# Least Absolute Deviation Utility for Trapezoidal Fuzzy Preference Relations

Lei He[a,b], Diego García-Zamora[c], Yuming Zhu[a*], Luis Martínez[b,d]

[a] *School of Management, Northwestern Polytechnical University, Xi'an, Shaanxi, 710072, P.R. China*
[b] *Computer Science Department, University of Jaén, Jaén, Jaén 23071, Spain*
[c] *Department of Mathematics, University of Jaén, Jaén, Jaén, 23071, Spain*
[d] *Andalusian Research Institute in Data Science and Computational Intelligence, Jaén, Jaén, 23071, Spain*

**Abstract:** Preference relations (PRs) are widely used to model expert judgments because they allow for eliciting the decision-makers' opinions from pairwise comparisons. Traditionally, PRs have been elicited using real numbers. However, in real-world decision-makers usually feel more comfortable using linguistic expressions closer to natural language. In this context, the purpose of this paper is to extend the classical idea of PR into the environment of Trapezoidal Fuzzy Numbers (TrFNs) by addressing several drawbacks in current research. On the one hand, existing fuzzy extensions for Fuzzy Preference Relations (FPRs) and Multiplicative Preference Relations (MPRs) assume that the notion of neutrality must be modeled by a crisp real number (i.e., 0.5 or 1), which fails to capture the subjective and diverse ways in which decision-makers may perceive indifference. On the other hand, from the methodological point of view, current research fails at providing a theoretical framework that unifies both FPRs and MPRs for fuzzy numbers and that simultaneously develops a coherent generalization of the involved classical notions (e.g., reciprocity, consistency) from the perspective of fuzzy arithmetic. Therefore, we begin by introducing the notion of neutral TrFN, a fuzzy number that can model the notion of indifference on a certain fuzzy scale. This generalization provides a flexible representation of neutrality, leading to a more realistic extension of the notion of PRs into the fuzzy environment. Building upon this, we propose the notion of Trapezoidal FPR (TrFPR), analyze its theoretical properties and study how reciprocity and consistency are extended using fuzzy arithmetic. Furthermore, we introduce the notion of Least Absolute Deviation (LAD) utility vector associated with a consistent TrFPR and develop an optimization model to derive it from an inconsistent TrFPR in order to enable the computation of priority values for TrFPRs. These ideas are also extended to MPRs in the fuzzy domain. Finally, the applicability of the proposed LAD utility method is demonstrated within a problem for evaluating land development projects, where comparison with fuzzy AHP illustrates the superiority of our proposal.

## 1. Introduction

Preference relations (PRs) are widely used to model experts' opinions in decision-making problems (Sahoo et al., 2023; Xu, 2014), as they allow experts to express their judgments through structured pairwise comparisons, which may provide more precise modeling to the human psyche than requiring direct evaluations of all alternatives individually (Jensen, 1984; Miller, 1956). Usually, decision-makers (DMs) express their pairwise comparison information in two formats: Fuzzy Preference Relations (FPRs) (Orlovsky, 1978) and Multiplicative Preference Relations (MPRs) (Saaty, 1990). The judgments in MPRs are quantified by positive ratios, and the judgments in FPRs are quantified by numeric values within the interval [0,1].

---

[*] Corresponding author
*E-mail addresses*: helei0131@mail.nwpu.edu.cn (Lei He), dgzamora@ujaen.es (Diego García-Zamora), zym1886@nwpu.edu.cn (Yuming Zhu), martin@ujaen.es (Luis Martínez).

To address the uncertainty and vagueness inherent in real-world decision-making processes (Pedrycz and Song, 2014), researchers have introduced fuzzy set theory (Zadeh, 1965) into the study of PRs and have subsequently proposed various extended models, such as interval-valued preference relations (Saaty and Vargas, 1987), hesitant fuzzy preference relations (Xia and Xu, 2013; Yuan Xu et al., 2022), and trapezoidal fuzzy preference relations (Meng and Chen, 2021). Compared with other forms of PRs, Trapezoidal Fuzzy Numbers (TrFNs) provide a flexible and intuitive approach to modeling uncertain preferences, and they can be easily derived from linguistic terms. Therefore, Trapezoidal Fuzzy Preference Relations (TrFPRs) (Gong et al., 2012) and Trapezoidal Multiplicative Preference Relations (TrMPRs) (Buckley, 1985) are particularly practical in complex decision-making environments (Meng et al., 2021; Wu et al., 2019; Zhou et al., 2018) and have garnered considerable attention in multi-criteria decision making (MCDM), group decision making (GDM), and risk analysis (Gong et al., 2012; Meng and Chen, 2021; Wan et al., 2023; Zhang et al., 2025).

Although TrFPRs and TrMPRs enhance the fuzzy expressiveness of PRs, existing studies generally extend their definitions from the classical FPRs and MPRs (Wu et al., 2019). A key limitation of these definitions is that they rely on a fixed, crisp neutral value (0.5 for FPRs and 1 for MPRs). However, this fixed representation of neutrality fails to capture the nuanced complexity of human judgments. In reality, DMs may perceive "indifference" as a fuzzy number rather than a precise point, and interpret uncertainty in various ways (García-Zamora et al., 2024; Jiang et al., 2022). For instance, when using methods such as fuzzy envelopes to translate linguistic terms to TrFNs, the neutral value is not necessarily fixed at 0.5 (Liu and Rodríguez, 2014). Therefore, this paper aims to move beyond this rigid assumption by considering the fuzziness of neutral preferences and proposing a more general definition of TrFPRs.

In addition to defining TrFPRs and analyzing their theoretical properties, another fundamental challenge in decision-making is to derive a priority vector from a PR to rank alternatives (Wang et al., 2007; Xu et al., 2016). For MPRs, prominent methods include the eigenvector method (Saaty, 2008) and the Best-Worst Method (BWM) (Rezaei, 2015). In addition, several researchers have developed optimization methods that derive priorities from FPRs by considering consistency (Chang et al., 2025; Wan et al., 2023). This approach typically employs optimization methods to find the consistent FPR that best approximates an inconsistent FPR (Wu et al., 2019). Following these ideas, we develop here an optimization-based method to obtain the priority vector from a TrFPR that enables the ranking of alternatives using fuzzy numbers. In this regard, we aim at tackling challenges such as (i) the proposal of a clear framework that relates multiplicative and additive fuzzy scales, as well as provides a differentiation between the notions of fuzzy weighting vector and fuzzy priority value (priority in Wang et al.(2007); weight in Wang(2015); priority weight in Xu et al.,(2016)); (ii) developing a practical framework based in fuzzy arithmetic and that enables the derivation of $n$ fuzzy priority values from a matrix with $n^2$ pairwise comparisons; (iii) integrating the notion of personalized neutral preference for TrFPRs in a rational way.

The main contributions of this study are as follows:
(1) We introduce the notion of personalized neutral TrFN to model the concept of indifference.
(2) We present a theoretical framework to extend classical FPRs and MPRs into the context of TrFNs, which is consistent with fuzzy philosophy and does not rely on a defuzzification mechanism that implies an information loss.
(3) We develop a method that allows mapping the proposed TrFPRs into fuzzy priority values, thereby reducing dimensionality and facilitating the interpretability of the results and the computations.
(4) We provide a unified framework to handle both additive and multiplicative TrFPRs, thus generalizing the derivation process across different TrFN scales.

As a consequence, our study provides a general framework for deriving utility values from TrFPRs, which not only considers a more general definition of TrFPRs but also integrates both additive and multiplicative preference models. In addition, this framework is not limited to specific representations such as hesitant fuzzy linguistic preference (Li et al., 2024), but rather allows any TrFN to be used as input in the TrFPR and also permits any symmetric TrFN to represent neutrality.

The remainder of this paper is organized as follows. Section 2 introduces the preliminaries, including basic concepts of fuzzy sets and PRs. Section 3 presents the definitions of the neutral TrFN and the TrFPR, then develops the LAD utility models for both TrFPRs and TrMPRs. Section 4 illustrates the applicability of our LAD utility by embedding it into the Fuzzy AHP and compares the proposed approach with the classical fuzzy weighting methods. Finally, Section 5 concludes the paper and discusses future research.

## 2. Preliminaries

In this section, we present the main concepts on fuzzy sets and preference relations that form the theoretical foundation of this study.

### 2.1. Fuzzy sets and fuzzy numbers

The fuzzy set theory (Zadeh, 1965) establishes the mathematical basis for our approach, extending the classical concept of set membership to allow partial membership degrees. It is worth noting that every subset $A \subseteq \mathbb{R}$ can be described by the indicator function $\chi_A : \mathbb{R} \to \{0,1\}$, defined as follows:

$$\chi_A(x) = \begin{cases} 1 \text{ if } x \in A, \\ 0 \text{ if } x \notin A. \end{cases}$$

However, real-world phenomena often exhibit gradual transitions that crisp boundaries cannot adequately capture. Fuzzy set theory generalizes this concept by allowing membership degrees in the interval [0,1].

**Definition 1** (*Fuzzy Sets* (Zadeh, 1965)). Let $U$ be the universe of discourse. A mapping function $A: U \to [0,1]$ is called a fuzzy set on $U$. The value $A(x)$ represents the degree of membership of $x$ in $A$.

Fuzzy numbers are a specific class of fuzzy sets that generalize the concept of real numbers.

**Definition 2** (*Fuzzy Number* (Klir and Yuan, 1995)). A fuzzy set $A : \mathbb{R} \to [0,1]$ is said to be a fuzzy number if it satisfies:

(1) $A$ is a normal fuzzy set, that is, there existing at least one $x \in \mathbb{R}$ such that $A(x)=1$,

(2) For any $\alpha \in (0,1]$, the $\alpha$-level set $A_\alpha = \{x \in \mathbb{R} : A(x) \geq \alpha\}$ is a non-empty bounded closed interval,

(3) The support of $A$, defined as $\text{supp}(A) = \text{Cl}(\{x \in \mathbb{R} : A(x) > 0\})$, is bounded.

In real-world decision-making, capturing the inherent vagueness and imprecision of human judgment is a significant challenge. DMs often find it more natural and intuitive to express their assessments using linguistic terms (e.g., "Good," "Fair," "Poor") rather than crisp numerical values (Zadeh, 1975). The fuzzy linguistic approach (Herrera and Martinez, 2000) provides a formal framework for modeling such information. A key idea in this approach is to represent the semantics of each linguistic term with a fuzzy number. Among various types of fuzzy numbers, TrFNs have gained particular prominence due to their exceptional balance between expressive power and computational simplicity (Zadeh, 1975; Delgado et al., 1998). Their structure is intuitive, and their operational simplicity facilitates the development of tractable models. Therefore, this study adopts TrFNs as the fundamental tool for modeling preferences. A TrFN is formally defined as follows.

**Definition 3** (*Trapezoidal Fuzzy Numbers* (Dubois and Prade, 1980)). A TrFN is a fuzzy number denoted as $\tilde{T} \equiv T(a,b,c,d)$ with $a \leq b \leq c \leq d$, whose membership function $\tilde{T}$ is defined as:

$$\tilde{T}(x) = \begin{cases} 0, & \text{if } x \leq a \text{ or } x \geq d \\ \frac{x-a}{b-a}, & \text{if } a < x < b \\ 1, & \text{if } b \leq x \leq c \\ \frac{d-x}{d-c}, & \text{if } c < x < d \end{cases}$$

For the sake of clarity, we use $\mathcal{T}(U)$ to denote the set of all TrFNs defined on a universe $U$. Specifically, let $\mathcal{T}([0,1])$, $\mathcal{T}(\mathbb{R})$, $\mathcal{T}(\mathbb{R}^+)$ and $\mathcal{T}(\mathbb{R}_0^+)$ represent the sets of TrFNs defined on the interval [0,1], the real line $\mathbb{R}$, the positive real line $\mathbb{R}^+ = \{x \in \mathbb{R} | x > 0\}$, and the non-negative real line $\mathbb{R}_0^+ = \{x \in \mathbb{R} | x \geq 0\}$, respectively. The corresponding $n$-fold Cartesian product will be denoted by $\mathcal{T}_n(U)$.

The similarity or dissimilarity between different fuzzy assessments can be quantified by comparing the distance between two TrFNs. One commonly used distance measure is the normalized Manhattan (L1) distance (Heilpern, 1997), which can be calculated between as follows:

$$d(\tilde{T}_1, \tilde{T}_2) = \frac{1}{4}(|a_1 - a_2| + |b_1 - b_2| + |c_1 - c_2| + |d_1 - d_2|).$$

where $\tilde{T}_1 \equiv T(a_1, b_1, c_1, d_1)$ and $\tilde{T}_2 \equiv T(a_2, b_2, c_2, d_2)$. Note that this form, based on absolute values rather than power operations, facilitates the linearization of the optimization models we intend to propose in the following section.

Furthermore, the arithmetic operations between $\tilde{T}_1$ and $\tilde{T}_2$ can be defined as follows (Klir and Yuan, 1995):

(1) $\tilde{T}_1 \oplus \tilde{T}_2 = (a_1 + a_2, b_1 + b_2, c_1 + c_2, d_1 + d_2)$,

(2) $\tilde{T}_1 \ominus \tilde{T}_2 = (a_1 - d_2, b_1 - c_2, c_1 - b_2, d_1 - a_2)$,

(3) $r\tilde{T}_1 = (ra_1, rb_1, rc_1, rd_1), r > 0$.

To facilitate the final decision-making, TrFNs can be compared using the method proposed by Abbasbandy and Hajjari (2009). This method translates the TrFN $\tilde{T} \equiv T(a,b,c,d)$ into a numerical value called magnitude, which is defined as follows:

$$Mag(\tilde{T}) = \frac{a + 5b + 5c + d}{12}.$$

Then, the ranking of $\tilde{T}_1$ and $\tilde{T}_2$ can be defined by the $Mag(\cdot)$ as follows:

i) $Mag(\tilde{T}_1) > Mag(\tilde{T}_2)$ iff $\tilde{T}_1 \succ \tilde{T}_2$,

ii) $Mag(\tilde{T}_1) < Mag(\tilde{T}_2)$ iff $\tilde{T}_1 \prec \tilde{T}_2$,

iii) $Mag(\tilde{T}_1) = Mag(\tilde{T}_2)$ iff $\tilde{T}_1 \sim \tilde{T}_2$.

### 2.2. Preference relations

In some decision-making problems, the subjective judgments of DMs are structured through pairwise comparisons between alternatives (Sahoo et al., 2023; Xu, 2014). This formal structure is known as a preference relation. The two most commonly adopted forms of preference relations are FPRs and MPRs, which provide different scales for representing DMs' preferences.

**Definition 4** (*Fuzzy preference relation* (*FPR*) (Orlovsky, 1978; Tanino, 1984)). Let $n \in \mathbb{N}$. A matrix $X = (x_{ij}) \in \mathcal{M}_{n \times n}([0,1])$ is called an FPR if it satisfies $x_{ij} + x_{ji} = 1, \forall i, j \in \{1, ..., n\}$. The set of all the FPRs of dimension $n$ is denoted as:

$$\mathbb{A}_n = \left\{ X \in \mathbb{M}_{n \times n}([0,1]) : x_{ij} + x_{ji} = 1 \ \forall i, j \in \{1, ..., n\} \right\}.$$

**Definition 5** (*Multiplicative preference relation* (*MPR*) (Saaty, 2008)). Let $n \in \mathbb{N}$. A matrix $Y = (y_{ij}) \in \mathcal{M}_{n \times n}(\mathbb{R}^+)$ is called an MPR if it satisfies $y_{ij} \cdot y_{ji} = 1, \forall i, j \in \{1, ..., n\}$. The preference intensities $y_{ij}$ are measured using a ratio scale $m \in \mathbb{N} \setminus \{1\}$, indicating that $y_{ij} \in [\frac{1}{m}, m]$, $\forall i, j = \{1, ..., n\}$. Usually, $m=9$, but it could be any other real number greater than one. The set of all the MPRs of dimension $n$ is denoted as:

$$\mathbb{M}_n = \{Y \in \mathcal{M}_{n \times n}(\mathbb{R}^+) : y_{ij} \cdot y_{ji} = 1, y_{ij} > 0 \ \forall i, j = \{1, ..., n\}\},$$

whereas the MPRs elicited in the scale $m$ will be denoted as:

$$\mathbb{M}_n^m = \{Y \in \mathcal{M}_{n \times n}(\mathbb{R}^+) : y_{ij} \cdot y_{ji} = 1, y_{ij} \in [\tfrac{1}{m}, m] \ \forall i, j = \{1, ..., n\}\}.$$

When considering preference relations, it is necessary to consider their consistency (Wan et al., 2023). In general, a pairwise comparison matrix is considered consistent when it does not include any contradictory evaluations (Saaty, 2008). Specifically, additive FPRs achieve consistency when the preference degrees satisfy an additive transitivity condition (Definition 6), whereas MPRs satisfy a multiplicative transitivity condition (Definition 7).

**Definition 6** (*Consistent FPR* (Herrera-Viedma et al., 2002)). An FPR $X \in \mathbb{A}_n$ is said to be consistent if it satisfies

$$x_{ij} + \tfrac{1}{2} = x_{ik} + x_{kj} \ \forall i, j, k = \{1, ..., n\}.$$

**Definition 7** (*Consistent MPR* (Saaty, 2008)). An MPR $Y \in \mathbb{M}_n^m$ is said to be consistent if it satisfies

$$y_{ij} = y_{ik} \cdot y_{kj} \ \forall i, j, k = \{1, ..., n\}.$$

Definitions 6 and 7 provide consistency conditions for FPR and MPR, although they are defined on different scales, they can be equivalent through a bijective mapping function as formalized in the following Proposition 1.

**Proposition 1** (García-Zamora and Martínez, 2025; Herrera-Viedma et al., 2002). Let $n \in \mathbb{N}$ and $m \in \mathbb{N} \setminus \{1\}$. The mapping $\Phi : \mathbb{A}_n \to \mathbb{M}_n^m$ defined as $\Phi(X) = (\phi(x_{ij})) = m^{2x_{ij}-1}$, $\forall X \in \mathbb{A}_n$, is a bijection whose inverse is given by the mapping $\Phi^{-1}: \mathbb{M}_n^m \to \mathbb{A}_n$ defined as $\Phi^{-1}(Y) = (\phi^{-1}(y_{ij})) = \frac{1}{2} + \frac{1}{2}\log_m y_{ij}$, $\forall Y \in \mathbb{M}_n^m$.

These mappings $\Phi$ and $\Phi^{-1}$ preserve the consistency property, which can be formalized as following proposition.

**Proposition 2** (García-Zamora and Martínez, 2025). An FPR $X \in \mathbb{A}_n$ is consistent if and only if the MPR $\Phi(X) \in \mathbb{M}_n^m$ is consistent.

This inherent equivalence allows us to derive the consistency condition for MPRs directly from the additive one.

**Example 1.** Let us illustrate the mappings $\Phi$ and $\Phi^{-1}$ with $m=9$. Let an additive preference $x_{ij} = 0.75$, then $y_{ij} = \phi(0.75) = 9^{2(0.75)-1} = 9^{0.5} = 3$. Correspondingly, the reciprocal $x_{ji} = 0.25$ is mapped to $y_{ji} = \phi(0.25) = 9^{2(0.25)-1} = 9^{-0.5} = \frac{1}{3}$. Note that reciprocity is preserved, as $y_{ij} \cdot y_{ji} = 3 \cdot \frac{1}{3} = 1$. Applying the inverse mapping $\Phi^{-1}$, $\phi^{-1}(y_{ij}) = \phi^{-1}(3) = \frac{1}{2} + \frac{1}{2}(\log_9 3) = 0.75 = x_{ij}$, confirming that $\Phi$ maps preferences between the two scales while maintaining their structural properties.

## 3. Least absolute deviations utility

Based on the definitions of FPRs and TrFNs, this study adopts TrFPRs to better capture the uncertainty and fuzziness inherent in human decision-making. Distinct from traditional approaches, this section proposes a more general definition of TrFPR by incorporating a neutral TrFN $\tilde{T}_0$ to represent DMs' subjective neutral preferences. Subsequently, we introduce the concept of LAD utility vector and develop a corresponding optimization model to derive it.

### 3.1. Trapezoidal fuzzy preference relations

To facilitate a more concise definition of TrFPR, we extend the notion of negation operator to TrFNs. In fuzzy set theory, negation functions serve as fundamental operators for modeling complement concepts. A function $N:[0,1] \to [0,1]$ is called a strong negation if it is decreasing and satisfies $N(N(x)) = x$ for all $x \in [0,1]$. In this paper, we use the standard negation $N_S(x) = 1-x$ (Zadeh, 1965), which serves as the canonical model for all strong negations (Trillas, 1979). Based on this function, we can define a corresponding standard negation operator for the set of TrFNs.

**Definition 8** (*Standard Negation for TrFNs*). Let $\tilde{T} \equiv T(a,b,c,d)$ be a TrFN defined on [0,1]. The standard negation of a TrFN $\tilde{T}$ is denoted as $\tilde{T}^\circ$ and defined by

$$\tilde{T}^\circ = T^\circ(a,b,c,d) = T(N_S(d), N_S(c), N_S(b), N_S(a)) = T(1-d, 1-c, 1-b, 1-a).$$

**Proposition 3** (*Involutive Property*). The standard negation operator $(\cdot)^\circ$ is involutive. For any $\tilde{T} \in \mathcal{T}([0,1])$, it holds that $(\tilde{T}^\circ)^\circ = \tilde{T}$.

**Proof.** Let $\tilde{T} \equiv T(a,b,c,d)$. The first application of the operator gives $T^\circ(a,b,c,d) = T(1-d, 1-c, 1-b, 1-a)$. Applying the operator again yields:

$$\begin{aligned}(\tilde{T}^\circ)^\circ &= (T(1-d, 1-c, 1-b, 1-a))^\circ \\ &= T(1-(1-a), 1-(1-b), 1-(1-c), 1-(1-d)) \\ &= T(a,b,c,d) = \tilde{T}.\end{aligned}$$

$\square$

**Example 2.** Let us consider a TrFN $\tilde{T} = (0.6, 0.7, 0.8, 0.9)$ representing a preference degree. According to Definition 8, we apply the $N_S(x) = 1-x$ to compute its standard negation $\tilde{T}^\circ$:

$$\tilde{T}^\circ = T(1-0.9,\ 1-0.8,\ 1-0.7,\ 1-0.6) = T(0.1,\ 0.2,\ 0.3,\ 0.4).$$

We can verify the involutive property by computing $(\tilde{T}°)° = T(1-0.4, 1-0.3, 1-0.2, 1-0.1) = T(0.6, 0.7, 0.8, 0.9) = \tilde{T}$.

In fuzzy set theory, a fundamental property of any continuous strong negation is the existence of a unique fixed point $e \in [0,1]$ such that $N(e) = e$ (Trillas, 1979). This fixed point represents the neutral value in the context of the negation. For the standard negation $N_S(x)$, its fixed point is $\frac{1}{2}$. Therefore, in the classical definition of FPR, the neutral reference point, representing "indifference," is fixed at the crisp value of $\frac{1}{2}$. As researchers extended these concepts to fuzzy number domains, the neutral preference was maintained as the crisp value $\frac{1}{2}$ (Meng et al., 2021). This rigid assumption overlooks the nuanced and often subjective nature of human judgment, where DMs may perceive and interpret uncertainty in different ways (García-Zamora et al., 2024; Jiang et al., 2022). A DM's sense of "perfect neutrality" might not be a single point but rather a fuzzy number, reflecting a zone of indifference. To address this limitation and build a more flexible framework, we introduce the concept of a neutral TrFN, which is grounded in the fundamental principle that a neutral element should be a fixed point under its corresponding negation operation. Therefore, based on the standard negation operator, $(\cdot)°$, we define a neutral TrFN $\tilde{T}_0$ to represent the DM's fuzzy perception of the neutral value.

**Definition 9** (*Neutral TrFN*). A TrFN $\tilde{T}_0 = T(T_0^a, T_0^b, T_0^c, T_0^d) \in \mathcal{T}([0,1])$ is called a neutral TrFN if it satisfies the fixed-point property with respect to the negation operator:

$$\tilde{T}_0° = \tilde{T}_0.$$

This condition ensures that the concept of neutrality represented by $\tilde{T}_0$ is invariant under the negation operation, which preserves the fixed-point property of strong negation. In the decision-making process, $\tilde{T}_0$ can be obtained, for example, by a fuzzy envelope (Liu and Rodríguez, 2014), or capturing the DM's fuzzy perception of neutrality by the deck of cards method interaction with the analyst (García-Zamora et al., 2024). Having established the notion of a neutral TrFN, it is essential to adopt a ranking method that is compatible with this new definition of neutrality. Specifically, the ranking method should preserve the property derived from the standard negation. To achieve this, we propose a generalized magnitude function, inspired by the work of Abbasbandy and Hajjari (2009).

**Definition 10** (*Magnitude of a TrFN*). For a TrFN $\tilde{T} \equiv T(a,b,c,d)$, its magnitude $Mag_w(\tilde{T})$ is defined as $Mag_w(\tilde{T}) = w_1(a+d) + w_2(b+c)$, where the weights $w_1, w_2 > 0$ satisfy $2(w_1 + w_2) = 1$.

The flexibility of $w_1$ and $w_2$ allows us to define a family of ranking functions. The crucial property of this magnitude function is its behavior under negation, as stated in the following proposition.

**Proposition 4.** Let us consider a TrFN $\tilde{T} \equiv T(a,b,c,d)$ with its standard negation $\tilde{T}°$ and a neutral TrFN $\tilde{T}_0 \in \mathcal{T}([0,1])$. Then $Mag_w(\tilde{T}°) = 1 - Mag_w(\tilde{T})$ and $Mag_w(\tilde{T}_0) = \frac{1}{2}$.

**Proof.** For a TrFN $\tilde{T} \equiv T(a,b,c,d)$, by Definition 8 we obtain

$$\begin{aligned}Mag_w(\tilde{T}°) &= Mag_w(T(1-d, 1-c, 1-b, 1-a)) = w_1((1-d)+(1-a)) + w_2((1-c)+(1-b)) \\ &= w_1(2-(a+d)) + w_2(2-(b+c)) \\ &= 2(w_1+w_2) - (w_1(a+d) + w_2(b+c)) \\ &= 1 - Mag_w(\tilde{T}).\end{aligned}$$

For the neutral TrFN $\tilde{T}_0$, the fixed-point property $\tilde{T}_0° = \tilde{T}_0$ implies: $Mag_w(\tilde{T}_0°) = 1 - Mag_w(\tilde{T}_0) = Mag_w(\tilde{T}_0)$. Thus, $Mag_w(\tilde{T}_0) = \frac{1}{2}$. □

This finding demonstrates that our proposed neutral TrFN, although defined as a TrFN with flexible and uncertain boundaries, still satisfies the fixed-point property of strong negation $Mag_w(\tilde{T}_0) = \frac{1}{2}$ and can thus be regarded as an effective measure of neutrality. Based on this, we propose the definition of TrFPR:

**Definition 11** (*Trapezoidal fuzzy preference relation (TrFPR)*). Let $n \in \mathbb{N}$ and let $\tilde{T}_0 \in \mathcal{T}([0,1])$ be a TrFN representing a fuzzy neutral preference value. For a preference matrix $\tilde{X} = (\tilde{x}_{ij}) \in \mathcal{M}_{n \times n}(\mathcal{T}([0,1])), \forall i, j = \{1,\ldots,n\}$, where $\tilde{x}_{ij} = (a_{ij}, b_{ij}, c_{ij}, d_{ij}) \in \mathcal{T}([0,1])$ is a TrFN, we will say that $\tilde{X}$ is a TrFPR if

(i) $\tilde{x}_{kk} = \tilde{T}_0 \ \forall k \in \{1,\ldots,n\}$,

(ii) $\tilde{x}_{ji} = \tilde{x}_{ij}^\circ \ \forall i, j \in \{1,\ldots,n\}, i \neq j$,

The set of TrFPRs of dimension $n$ with neutral TrFN $\tilde{T}_0$ is denoted as $\mathbb{A}_n(\tilde{T}_0)$.

**Example 3.** A TrFPR $\tilde{X}$ for $n = 3$ alternatives with the neutral TrFN $\tilde{T}_0 = T(0.4, 0.5, 0.5, 0.6)$ could be constructed as follows:

$$\tilde{X} = \begin{pmatrix} T(0.4,0.5,0.5,0.6) & T(0.6,0.7,0.7,0.8) & T(0.6,0.7,0.8,0.9) \\ T(0.2,0.3,0.3,0.4) & T(0.4,0.5,0.5,0.6) & T(0.5,0.6,0.7,0.8) \\ T(0.1,0.2,0.3,0.4) & T(0.2,0.3,0.4,0.5) & T(0.4,0.5,0.5,0.6) \end{pmatrix}.$$

This matrix satisfies: i) $\tilde{x}_{11} = \tilde{x}_{22} = \tilde{x}_{33} = T(0.4, 0.5, 0.5, 0.6) = \tilde{T}_0$. ii) $\tilde{x}_{12}^\circ = (T(0.6, 0.7, 0.7, 0.8))^\circ = T(1-0.8, 1-0.7, 1-0.7, 1-0.6) = T(0.2, 0.3, 0.3, 0.4) = \tilde{x}_{21}$. Similarly, $\tilde{x}_{13}^\circ = \tilde{x}_{31}$ and $\tilde{x}_{23}^\circ = \tilde{x}_{32}$.

**Definition 12** (*Consistent TrFPR*). A TrFPR $\tilde{X} \in \mathbb{A}_n(\tilde{T}_0)$ is said to be consistent if it satisfies

$$\tilde{x}_{ij} \oplus \tilde{T}_0 = \tilde{x}_{ik} \oplus \tilde{x}_{kj} \ \forall i, j, k \in \{1,\ldots,n\},$$

where $\tilde{T}_0 \in \mathcal{T}([0,1])$ represents the fuzzy neutral preference and $\oplus$ denotes the addition operation of TrFN (according to the operational rules defined in Definition 3).

If the fuzzy neutral preference degenerates into a crisp value (e.g., $\tilde{T}_0 = T(\frac{1}{2}, \frac{1}{2}, \frac{1}{2}, \frac{1}{2})$), and all input preference values $\tilde{x}_{ij}$ are crisp numbers, then both the definition of TrFPR and consistent TrFPR reduce to the classical definitions of FPR (Orlovsky, 1978) and consistent FPR (Herrera-Viedma et al., 2002).

**Example 4.** Continuing with the TrFPR $\tilde{X}$ and the neutral TrFN $\tilde{T}_0$ from Example 3, note that a consistent TrFPR should satisfy $\tilde{x}_{ij} \oplus \tilde{T}_0 = \tilde{x}_{ik} \oplus \tilde{x}_{kj}$ for all $i, j, k \in \{1,\ldots,n\}$. For the case $i=1, j=2, k=3$,

$$\tilde{x}_{12} \oplus \tilde{T}_0 = T(0.6, 0.7, 0.7, 0.8) \oplus T(0.4, 0.5, 0.5, 0.6) = T(1.0, 1.2, 1.2, 1.4),$$
$$\tilde{x}_{13} \oplus \tilde{x}_{32} = T(0.6, 0.7, 0.8, 0.9) \oplus T(0.2, 0.3, 0.4, 0.5) = T(0.8, 1.0, 1.2, 1.4).$$

It follows that $\tilde{x}_{ij} \oplus \tilde{T}_0 \neq \tilde{x}_{ik} \oplus \tilde{x}_{kj}$, thus we conclude that the given TrFPR matrix $\tilde{X}$ is not consistent. Such inconsistency is a common characteristic of subjectively elicited preference information, highlighting the necessity of a robust methodology to derive a rational and consistent utility vector from imperfect or uncertain data.

### 3.2. Least absolute deviations utility for TrFPRs

While a TrFPR matrix $\tilde{X}$ comprehensively captures pairwise preferences, its $n \times n$ structure is often difficult to interpret directly, making it challenging to determine an overall ranking of alternatives. To address this, our goal is to distill this complex relational data into a more concise and intuitive fuzzy utility vector, $\tilde{u} = (\tilde{u}_1, \ldots, \tilde{u}_n) \in \mathcal{T}_n([0,1])$. The fuzzy utility vector not only offers superior interpretability by providing a clear ranking but also enhances computational efficiency by reducing the dimensionality from $n^2$ to $n$, which is particularly beneficial for handling large-scale data (García-Zamora and Martínez, 2025). Therefore, this study focuses on developing a method to derive the fuzzy utility vector from a TrFPR.

The theoretical foundation for deriving the fuzzy utility vector lies in the direct correspondence between a utility vector and a consistent TrFPR. This relationship is formally established in the following theorem.

**Theorem 1.** Let $n \in \mathbb{N}$. A TrFPR $\tilde{X} \in \mathbb{A}_n(\tilde{T}_0)$ is consistent if and only if there exist a trapezoidal fuzzy utility vector $\tilde{u} \in \mathcal{T}_n([0,1])$ such that
$$\tilde{x}_{ij} \oplus \tilde{T}_0 = \tilde{u}_i \oplus \tilde{u}_j^\circ \quad \forall i,j = \{1,...,n\}.$$

**Proof.** i) To prove the sufficient condition, if the TrFPR $\tilde{X} \in \mathbb{A}_n(\tilde{T}_0)$ is consistent, let us pick $\tilde{u}_i := \tilde{x}_{ik}$ and $\tilde{u}_j := \tilde{x}_{jk}$ for some $k = \{1,...,n\}$. Then, $\tilde{u}_i \oplus \tilde{u}_j^\circ = \tilde{x}_{ik} \oplus \tilde{x}_{jk}^\circ = \tilde{x}_{ik} \oplus \tilde{x}_{kj} = \tilde{x}_{ij} \oplus \tilde{T}_0$. Therefore, there exists a trapezoidal fuzzy utility vector $\tilde{u} \in \mathcal{T}_n([0,1])$ such that $\tilde{u}_i \oplus \tilde{u}_j^\circ = \tilde{x}_{ij} \oplus \tilde{T}_0$.

ii) To prove the necessary condition, let us assume that there exists a trapezoidal fuzzy utility vector $\tilde{u} \in \mathcal{T}_n([0,1])$ such that $\tilde{x}_{ij} \oplus \tilde{T}_0 = \tilde{u}_i \oplus \tilde{u}_j^\circ, \forall i,j = \{1,...,n\}$. In such case, $\tilde{x}_{ik} \oplus \tilde{x}_{kj} \oplus 2\tilde{T}_0 = \tilde{x}_{ik} \oplus \tilde{T}_0 \oplus \tilde{x}_{kj} \oplus \tilde{T}_0 = \tilde{u}_i \oplus \tilde{u}_k^\circ \oplus \tilde{u}_k \oplus \tilde{u}_j^\circ = \tilde{u}_i \oplus \tilde{u}_j^\circ \oplus 2\tilde{T}_0 = \tilde{x}_{ij} \oplus 3\tilde{T}_0$. Thus, $\tilde{x}_{ik} \oplus \tilde{x}_{kj} = \tilde{x}_{ij} \oplus \tilde{T}_0$, indicating that the TrFPR $\tilde{X} \in \mathbb{A}_n(\tilde{T}_0)$ is consistent. □

Theorem 1 establishes the relationship between a consistent TrFPR $\tilde{X}$ and its utility vector $\tilde{u}$. Since the primary goal of deriving a utility vector is to facilitate the ranking of alternatives, it is essential to ensure that this correspondence preserves the preference ordering. The following corollary verifies that the ranking of alternatives obtained from the utility vector $\tilde{u}$ is same with the pairwise preferences in $\tilde{X}$.

**Corollary 1.** Let $\tilde{X} \in \mathbb{A}_n(\tilde{T}_0)$ be a consistent TrFPR and $\tilde{u} = (\tilde{u}_1,...,\tilde{u}_n) \in \mathcal{T}_n([0,1])$ is one of the trapezoidal fuzzy utility vectors given by Theorem 1. Then for any $i,j = \{1,...,n\}$:

i) $Mag_w(\tilde{x}_{ij}) > Mag_w(\tilde{T}_0)$ iff $Mag_w(\tilde{u}_i) > Mag_w(\tilde{u}_j)$,

ii) $Mag_w(\tilde{x}_{ij}) = Mag_w(\tilde{T}_0)$ iff $Mag_w(\tilde{u}_i) = Mag_w(\tilde{u}_j)$,

iii) $Mag_w(\tilde{x}_{ij}) < Mag_w(\tilde{T}_0)$ iff $Mag_w(\tilde{u}_i) < Mag_w(\tilde{u}_j)$.

**Proof.** For a consistent TrFPR $\tilde{X} = (x_{ij})$ such that $\tilde{x}_{ij} \oplus \tilde{T}_0 = \tilde{u}_i \oplus \tilde{u}_j^\circ, \forall i,j \in \{1,...,n\}$, it holds that the $Mag_w(\tilde{x}_{ij}) + Mag_w(\tilde{T}_0) = Mag_w(\tilde{u}_i) + Mag_w(\tilde{u}_j^\circ)$. Therefore, by Proposition 4, it follows that:

$$Mag_w(\tilde{x}_{ij}) + Mag_w(\tilde{T}_0) = Mag_w(\tilde{u}_i) + (1 - Mag_w(\tilde{u}_j)) \forall i,j = \{1,...,n\}$$
$$\Leftrightarrow Mag_w(\tilde{x}_{ij}) + (1 - Mag_w(\tilde{T}_0)) = Mag_w(\tilde{u}_i) + (1 - Mag_w(\tilde{u}_j)) \forall i,j = \{1,...,n\}$$
$$\Leftrightarrow Mag_w(\tilde{x}_{ij}) - Mag_w(\tilde{T}_0) = Mag_w(\tilde{u}_i) - Mag_w(\tilde{u}_j) \forall i,j = \{1,...,n\}.$$

Then, we can prove the three relationships: i) $Mag_w(\tilde{x}_{ij}) > Mag_w(\tilde{T}_0)$ if and only if $Mag_w(\tilde{u}_i) > Mag_w(\tilde{u}_j)$, ii) $Mag_w(\tilde{x}_{ij}) = Mag_w(\tilde{T}_0)$ if and only if $Mag_w(\tilde{u}_i) = Mag_w(\tilde{u}_j)$, iii) $Mag_w(\tilde{x}_{ij}) < Mag_w(\tilde{T}_0)$ if and only if $Mag_w(\tilde{u}_i) < Mag_w(\tilde{u}_j)$. □

**Example 5.** Let us consider Example 3 with $n=2$ and $\tilde{T}_0 = T(0.4, 0.5, 0.5, 0.6)$. To illustrate a fully consistent TrFPR, we presuppose an underlying utility vector $\tilde{u} = (\tilde{u}_1, \tilde{u}_2, \tilde{u}_3)$, where $\tilde{u}_1 = T(0.3, 0.3, 0.3, 0.5)$, $\tilde{u}_2 = T(0.1, 0.1, 0.1, 0.3)$, and $\tilde{u}_3 = T(0.0, 0.0, 0.0, 0.2)$. This utility vector is, in fact, the optimal solution found for the inconsistent matrix in Example 6. Then, we can construct the fully consistent TrFPR $\tilde{X}^*$ by Theorem 1 as follows:

$$\tilde{X}^* = \begin{pmatrix} T(0.4,0.5,0.5,0.6) & T(0.6,0.7,0.7,0.8) & T(0.7,0.8,0.8,0.9) \\ T(0.2,0.3,0.3,0.4) & T(0.4,0.5,0.5,0.6) & T(0.5,0.6,0.6,0.7) \\ T(0.1,0.2,0.2,0.3) & T(0.3,0.4,0.4,0.5) & T(0.4,0.5,0.5,0.6) \end{pmatrix}.$$

For the case $i=1, j=2, k=3$,

$$\tilde{x}_{12} \oplus \tilde{T}_0 = T(0.6, 0.7, 0.7, 0.8) \oplus T(0.4, 0.5, 0.5, 0.6) = T(1.0, 1.2, 1.2, 1.4),$$
$$\tilde{x}_{13} \oplus \tilde{x}_{32} = T(0.7, 0.8, 0.8, 0.9) \oplus T(0.3, 0.4, 0.4, 0.5) = T(1.0, 1.2, 1.2, 1.4).$$

It follows that $\tilde{x}_{ij} \oplus \tilde{T}_0 = \tilde{x}_{ik} \oplus \tilde{x}_{kj}$, thus the TrFPR $\tilde{X}^*$ constructed from the utility vector $\tilde{u}$ is consistent.

Theorem 1 and Example 5 have shown that, given a consistent TrFPR $\tilde{X} \in \mathbb{A}_n(\tilde{T}_0)$, we can equivalently represent the shifted preference information, $\tilde{x}_{ij} \oplus \tilde{T}_0$, through a utility vector $\tilde{u} \in \mathcal{T}_n([0,1])$. Conversely, a utility vector $\tilde{u} \in \mathcal{T}_n([0,1])$ can be used to construct a consistent TrFPR. However, the preference matrices elicited from DMs are usually inconsistent in practice due to cognitive limitations and environmental uncertainty (Saaty, 2008; Wan et al., 2023). Therefore, it necessitates an alternative approach to derive the most representative utility vector from an inconsistent TrFPR.

To address it, we propose to derive a fuzzy utility vector $\tilde{u} = (\tilde{u}_1, ..., \tilde{u}_n) \in \mathcal{T}_n([0,1])$ that minimizes the distance between the shifted original TrFPR $\tilde{x}_{ij} \oplus \tilde{T}_0$ and the shifted preference relation reconstructed from the utility vector $\tilde{u}$ by $\tilde{u}_i \oplus \tilde{u}_j^\circ$. For this purpose, we adopt the normalized Manhattan (L1) distance for TrFNs, as defined in Definition 3, to compute the sum of absolute deviations based on the principle of Least Absolute Deviations (LAD). The LAD approach preserves trapezoidal semantics through direct parameter comparison, and its absolute-value structure maintains linear tractability, enabling efficient reformulation as a linear programming model. It is worth noting that $d(\tilde{x}_{ij} \oplus \tilde{T}_0, \tilde{u}_i \oplus \tilde{u}_j^\circ)$ involves the use of the standard negation operator $(\cdot)^\circ$ within $\tilde{u}_j^\circ$, which is formally defined on $\tilde{u} \in \mathcal{T}_n([0,1])$. However, in practice, a highly inconsistent TrFPR may yield an optimal utility vector with components outside [0,1] interval. To be able to account for such cases, we will extend the standard negation operator $(\cdot)^\circ$ to all TrFNs, $\mathcal{T}(\mathbb{R})$, by using the same definition. This is a reasonable extension of the notation, justified by the operation $\tilde{u}_j^\circ$ maintains its essential property as a reflection around the fixed point $\frac{1}{2}$ across the real line. With this understanding, we formulate the following baseline optimization model to derive the LAD utility vector for TrFPR.

For a given TrFPR $\tilde{X} = (\tilde{x}_{ij}) \in \mathbb{A}_n(\tilde{T}_0)$. A fuzzy utility vector $\tilde{u} = (\tilde{u}_1, ..., \tilde{u}_n) \in \mathcal{T}_n(\mathbb{R})$ with $\tilde{u}_k = T(u_k^a, u_k^b, u_k^c, u_k^d)$ can be obtained by solving the optimization model **P₀**.

$$\min_{\tilde{u}_1, ..., \tilde{u}_n \in \mathcal{T}_n(\mathbb{R})} \sum_{i=1}^n \sum_{j=1}^n d(\tilde{x}_{ij} \oplus \tilde{T}_0, \tilde{u}_i \oplus \tilde{u}_j^\circ)$$
$$s.t. \{u_k^a \leq u_k^b \leq u_k^c \leq u_k^d \quad \forall k = \{1, ..., n\}\}$$

(**P₀**)

While **P₀** provides a theoretical baseline, it has a practical limitation. For example, a fuzzy utility $\tilde{u}_3 = T(-0.1, -0.1, -0.1, 0.1)$ is derived from **P₀** using the TrFPR in Example 3 contains negative components. However, a utility, as a measure of preference intensity, should be non-negative, then we introduce a non-negativity constraint and propose model **P**.

$$\min_{\tilde{u}_1, ..., \tilde{u}_n \in \mathcal{T}_n(\mathbb{R}_0^+)} \sum_{i=1}^n \sum_{j=1}^n d(\tilde{x}_{ij} \oplus \tilde{T}_0, \tilde{u}_i \oplus \tilde{u}_j^\circ)$$
$$s.t. \{0 \leq u_k^a \leq u_k^b \leq u_k^c \leq u_k^d \quad \forall k = \{1, ..., n\}\}$$

(**P**)

The introduction of the non-negativity constraint does not affect the optimal objective value of model **P₀**, since model **P** can be regarded as identifying a non-negative solution within the optimal solution set of **P₀**. This property is formally established in Proposition 5.

**Proposition 5.** Let **P₀** denotes the base LAD model, and **P** denotes the non-negative LAD model. If $\tilde{u}^* = (\tilde{u}_1^*, ..., \tilde{u}_n^*)$ is an optimal solution of **P₀**, then for any $\delta \in \mathbb{R}_0^+$, the utility vector $\tilde{u}'(\delta) = (\tilde{u}_1^* \oplus \tilde{\delta}, ..., \tilde{u}_n^* \oplus \tilde{\delta})$ is also an optimal solution of **P₀**. In particular, if $\delta_0 = -\min\{u_1^{a*}, u_2^{*a}..., u_n^{*a}\}$, then $\tilde{u}'(\delta_0) = (\tilde{u}_1^* \oplus \tilde{\delta}_0, ..., \tilde{u}_n^* \oplus \tilde{\delta}_0)$ is an optimal solution for both **P₀** and **P**.

**Proof.** Let $\tilde{u}^* = (\tilde{u}_1^*, ..., \tilde{u}_n^*)$ be an optimal solution for **P₀** with the optimal objective value given by $\frac{1}{4}\sum_{i=1}^{n}\sum_{j=1}^{n} d\left(\tilde{x}_{ij} \oplus \tilde{T}_0, \tilde{u}_i^* \oplus (\tilde{u}_j^*)^\circ\right)$. Let us construct a new solution $\tilde{u}'(\delta) = (\tilde{u}_1^* \oplus \tilde{\delta}, ..., \tilde{u}_n^* \oplus \tilde{\delta})$, then the objective value for $\tilde{u}'(\delta)$ can be expressed as $\frac{1}{4}\sum_{i=1}^{n}\sum_{j=1}^{n} d\left(\tilde{x}_{ij} \oplus \tilde{T}_0, (\tilde{u}_i^* \oplus \tilde{\delta}) \oplus (\tilde{u}_j^* \oplus \tilde{\delta})^\circ\right)$. Applying to the Definition 8, it holds that

$$(\tilde{u}_j^* \oplus \tilde{\delta})^\circ = T(1-(u_j^{d*}+\delta), 1-(u_j^{c*}+\delta), 1-(u_j^{b*}+\delta), 1-(u_j^{a*}+\delta))$$
$$= T((1-u_j^{d*})-\delta, (1-u_j^{c*})-\delta, (1-u_j^{b*})-\delta, (1-u_j^{a*})-\delta))$$
$$= (\tilde{u}_j^*)^\circ \ominus \tilde{\delta}.$$

Thus, $\frac{1}{4}\sum_{i=1}^{n}\sum_{j=1}^{n}\left(\tilde{x}_{ij} \oplus \tilde{T}_0, (\tilde{u}_i^* \oplus \tilde{\delta}) \oplus ((\tilde{u}_j^*)^\circ \ominus \tilde{\delta})\right) = \frac{1}{4}\sum_{i=1}^{n}\sum_{j=1}^{n} d\left(\tilde{x}_{ij} \oplus \tilde{T}_0, \tilde{u}_i^* \oplus (\tilde{u}_j^*)^\circ\right)$, it indicates the solution $\tilde{u}'(\delta) = (\tilde{u}_1^* \oplus \tilde{\delta}, ..., \tilde{u}_n^* \oplus \tilde{\delta})$ is also an optimal solution of **P₀**. Now let us consider $\delta_0 = -\min\{u_1^{a*}, u_2^{*a}..., u_n^{*a}\}$. We have $u_i^{a'} = u_i^{a*} + \delta_0 = u_i^{a*} - \min\{u_1^{a*}, u_2^{a*}..., u_n^{a*}\}$, which ensures that $u_i^{a'} \geq 0$ and then $\tilde{u}'(\delta_0) \geq 0$. Thus, $\tilde{u}'(\delta_0)$ is a feasible solution of **P**. Since the feasible region of **P** is a subset of the feasible region of **P₀**, the optimal value of **P** cannot be better than that of **P₀**. Therefore, $\tilde{u}'(\delta_0)$ is an optimal solution for both **P₀** and **P**. □

Proposition 5 guarantees that **P** derives a meaningful, non-negative utility vector without any loss of optimality. Building upon model **P**, we formally define the LAD utility by extending the constraints to the [0,1] interval for enhanced interpretability and consistency.

**Definition 13** (*LAD utility for TrFPRs*). Let $\tilde{X} = (\tilde{x}_{ij}) \in \mathbb{A}_n(\tilde{T}_0)$ be a TrFPR. An LAD utility vector for $\tilde{X}$ is a trapezoidal fuzzy vector $\tilde{u} = (\tilde{u}_1, \tilde{u}_2, ..., \tilde{u}_n) \in \mathcal{T}_n([0,1])$ that minimizes the total absolute deviation between the shifted original TrFPR $\tilde{x}_{ij} \oplus \tilde{T}_0$ and that generated by the utility vector, i.e.,

$$\min_{\tilde{u}_1,...,\tilde{u}_n \in \mathcal{T}_n([0,1])} \sum_{i=1}^{n}\sum_{j=1}^{n} d\left(\tilde{x}_{ij} \oplus \tilde{T}_0, \tilde{u}_i \oplus \tilde{u}_j^\circ\right)$$

subject to:

$$0 \leq u_k^a \leq u_k^b \leq u_k^c \leq u_k^d \leq 1 \ \forall k = \{1, ..., n\},$$

where $d(\cdot, \cdot)$ is the normalized Manhattan (L1) distance for TrFNs (Definition 3), $\tilde{u}_j^\circ$ is the negation of $\tilde{u}_j$, and $\tilde{T}_0$ is the reference TrFN representing the neutral preference.

In Definition 13, the upper-bound constraint $u_k^d \leq 1$ is explicitly imposed primarily to ensure the practical interpretability and formal consistency of the results. This constraint restricts utility values to a meaningful and standardized [0,1] interval, thereby preventing the model from yielding potentially large, impractical values. If this constraint were relaxed, highly inconsistent preferences in $\tilde{X}$ could lead to optimal utility values exceeding 1. This phenomenon serves as a indicator of high inconsistency, indicating that the DM may be advised to revise their judgments.

**Example 6.** To illustrate the application of Definition 13, we consider the TrFPR $\tilde{X}$ from Example 3, with the neutral preference $\tilde{T}_0 = T(0.4, 0.5, 0.5, 0.6)$. By solving the LAD model **P** using linear programming **L-P**, we obtain the optimal utility vector: $\tilde{u}^* = (T(0.3, 0.3, 0.3, 0.5), T(0.1, 0.1, 0.1, 0.3), T(0, 0, 0, 0.2))$ with an objective value of 0.2. Using the magnitude function $Mag(\cdot)$, the utility vector yields the ranking $A_1 \succ A_2 \succ A_3$.

It is worth noting that in real-world decision-making problems such as resource allocation, the total utility vector is often constrained by the available total resources. Therefore, based on the model **P**, a total utility constraint can be introduced to reflect this limitation, where the total utility can be represented by a TrFN $\tilde{\sigma} \in \mathcal{T}(\mathbb{R}^+)$. Then a utility vector $\tilde{u} = (\tilde{u}_1, \tilde{u}_2, ..., \tilde{u}_n) \in \mathcal{T}_n([0,1])$ can be obtained by solving the optimization model $\mathbf{P}_\sigma$:

$$\min_{\tilde{u}_1,...,\tilde{u}_n \in \mathcal{T}_n(\mathbb{R}_0^+)} \sum_{i=1}^n \sum_{j=1}^n d\left(\tilde{x}_{ij} \oplus \tilde{T}_0, \tilde{u}_i \oplus \tilde{u}_j^\circ\right)$$

$$s.t. \begin{cases} \bigoplus_{k=1}^n \tilde{u}_k = \tilde{\sigma} \\ 0 \leq u_k^a \leq u_k^b \leq u_k^c \leq u_k^d \quad \forall k = \{1,...,n\} \end{cases} \quad (\mathbf{P}_\sigma)$$

Due to the complexity of the fuzzy constraints and the objective structure of these models ($\mathbf{P_0}$, **P**, and $\mathbf{P}_\sigma$), it is not easy to find closed-form analytical solutions, and they have to be solved using numerical optimization methods. Based on Definition 3 of the distance between two TrFNs, we can express the objective function in these models more explicitly. Let the original TrFPR be defined as $\tilde{x}_{ij} = (a_{ij}, b_{ij}, c_{ij}, d_{ij})$ and let the utility be $\tilde{u}_i = (u_i^a, u_i^b, u_i^c, u_i^d)$. Then, the reconstructed preference relation can be denoted as: $\hat{u}_{ij} = \tilde{u}_i \oplus \tilde{u}_j^\circ$, which gives the following component-wise expressions: $\hat{u}_{ij}^a = u_i^a + 1 - u_j^d, \hat{u}_{ij}^b = u_i^b + 1 - u_j^c, \hat{u}_{ij}^c = u_i^c + 1 - u_j^b, \hat{u}_{ij}^d = u_i^d + 1 - u_j^a$. However, the objective function $\min \frac{1}{4} \sum_{i=1}^n \sum_{j=1}^n \left(|a_{ij} + T_0^a - \hat{u}_{ij}^a| + |b_{ij} + T_0^b - \hat{u}_{ij}^b| + |c_{ij} + T_0^c - \hat{u}_{ij}^c| + |d_{ij} + T_0^d - \hat{u}_{ij}^d|\right)$, involves absolute value terms, and is a non-linear model. To enhance computational tractability, we transform it into an equivalent linear program as stated in the following theorem.

**Theorem 2.** The non-linear optimization model **P** can be equivalently transformed into a linear programming model **L-P**.

$$\min_{\tilde{u}_1,...,\tilde{u}_n \in \mathcal{T}_n(\mathbb{R}_0^+)} \frac{1}{4} \sum_{i=1}^n \sum_{j=1}^n (v_{ij}^a + v_{ij}^b + v_{ij}^c + v_{ij}^d)$$

$$s.t. \begin{cases} v_{ij}^a \geq a_{ij} + T_0^a - \hat{u}_{ij}^a \\ v_{ij}^a \geq \hat{u}_{ij}^a - a_{ij} - T_0^a \\ v_{ij}^b \geq b_{ij} + T_0^b - \hat{u}_{ij}^b \\ v_{ij}^b \geq \hat{u}_{ij}^b - b_{ij} - T_0^b \\ v_{ij}^c \geq c_{ij} + T_0^c - \hat{u}_{ij}^c \\ v_{ij}^c \geq \hat{u}_{ij}^c - c_{ij} - T_0^c \\ v_{ij}^d \geq d_{ij} + T_0^d - \hat{u}_{ij}^d \\ v_{ij}^d \geq \hat{u}_{ij}^d - d_{ij} - T_0^d \\ 0 \leq u_k^a \leq u_k^b \leq u_k^c \leq u_k^d \quad \forall k = \{1,...,n\} \end{cases} \quad (\mathbf{L\text{-}P})$$

where $\hat{u}_{ij}^a = u_i^a + 1 - u_j^d, \hat{u}_{ij}^b = u_i^b + 1 - u_j^c, \hat{u}_{ij}^c = u_i^c + 1 - u_j^b, \hat{u}_{ij}^d = u_i^d + 1 - u_j^a$, and $v_{ij}^a, v_{ij}^b, v_{ij}^c, v_{ij}^d$ are auxiliary variables representing the absolute deviations for each component $(a,b,c,d)$ of the TrFNs $\tilde{x}_{ij} = (a_{ij}, b_{ij}, c_{ij}, d_{ij})$.

**Proof.** The core of this transformation lies in reformulating the non-linear absolute value term within the objective function of model **P** into a set of linear constraints. By substituting $d(\cdot,\cdot)$ into the objective function of **P**, we obtain:

$$\min_{\tilde{u}_1,...,\tilde{u}_n \in \mathcal{T}_n([0,1])} \sum_{i=1}^n \sum_{j=1}^n d\left(\tilde{x}_{ij} \oplus \tilde{T}_0, \tilde{u}_i \oplus \tilde{u}_j^\circ\right) =$$

$$\min_{\tilde{u}_1,...,\tilde{u}_n \in \mathcal{T}_n([0,1])} \frac{1}{4} \sum_{i=1}^n \sum_{j=1}^n \left(|a_{ij} + T_0^a - \hat{u}_{ij}^a| + |b_{ij} + T_0^b - \hat{u}_{ij}^b| + |c_{ij} + T_0^c - \hat{u}_{ij}^c| + |d_j + T_0^d - \hat{u}_{ij}^d|\right).$$

where $\hat{u}_{ij}^a = u_i^a + 1 - u_j^d, \hat{u}_{ij}^b = u_i^b + 1 - u_j^c, \hat{u}_{ij}^c = u_i^c + 1 - u_j^b, \hat{u}_{ij}^d = u_i^d + 1 - u_j^a$. To linearize each term $|\alpha_{ij} + T_0^\alpha - \hat{u}_{ij}^\alpha|$ for $\alpha \in \{a,b,c,d\}$, we introduce auxiliary variables $v_{ij}^\alpha \in \mathbb{R}_0^+$ constrained by:

$$\begin{cases} v_{ij}^\alpha \geq \alpha_{ij} + T_0^\alpha - \hat{u}_{ij}^\alpha \\ v_{ij}^\alpha \geq \hat{u}_{ij}^\alpha - \alpha_{ij} - T_0^\alpha \end{cases} \forall i,j \in \{1,...,n\}, \alpha \in \{a,b,c,d\}.$$

These constraints ensure that $v_{ij}^\alpha \geq |\alpha_{ij} + T_0^\alpha - \hat{u}_{ij}^\alpha|$ for all $i,j \in \{1,...,n\}$ and $\alpha \in \{a,b,c,d\}$. In **L-P**, the objective function minimizes $\frac{1}{4}\sum_{i=1}^n \sum_{j=1}^n (v_{ij}^a + v_{ij}^b + v_{ij}^c + v_{ij}^d)$ subject to these constraints, which guarantees equality at optimum: $v_{ij}^{\alpha*} = |\alpha_{ij} + T_0^{\alpha*} - \hat{u}_{ij}^{\alpha*}|$. Consequently, the optimal objective value of **L-P** can be expressed as

$$\tfrac{1}{4}\sum_{i=1}^n \sum_{j=1}^n \left( \left|a_{ij} + T_0^a - \hat{u}_{ij}^a\right| + \left|b_{ij} + T_0^b - \hat{u}_{ij}^b\right| + \left|c_{ij} + T_0^c - \hat{u}_{ij}^c\right| + \left|d_{ij} + T_0^d - \hat{u}_{ij}^d\right| \right),$$

which is equal to that of **P**. Since the constraints on the utility vector $\tilde{u}$ are identical in both models, and their objective functions yield the same optimal value, the non-linear model **P** and the linear model **L-P** are equivalent. Solving the linear model **L-P** will yield the same optimal utility vector $\tilde{u}^*$ as solving the non-linear model **P**. □

Similarly, the linearized models **L-P₀** and **L-P_σ** can also be equivalently transformed from models **P₀** and **P_σ** as following.

$$\min_{\tilde{u}_1,...,\tilde{u}_n \in \mathcal{T}_n(\mathbb{R})} \tfrac{1}{4}\sum_{i=1}^n \sum_{j=1}^n (v_{ij}^a + v_{ij}^b + v_{ij}^c + v_{ij}^d)$$

$$s.t. \begin{cases} v_{ij}^a \geq a_{ij} + T_0^a - \hat{u}_{ij}^a \\ v_{ij}^a \geq \hat{u}_{ij}^a - a_{ij} - T_0^a \\ v_{ij}^b \geq b_{ij} + T_0^b - \hat{u}_{ij}^b \\ v_{ij}^b \geq \hat{u}_{ij}^b - b_{ij} - T_0^b \\ v_{ij}^c \geq c_{ij} + T_0^c - \hat{u}_{ij}^c \\ v_{ij}^c \geq \hat{u}_{ij}^c - c_{ij} - T_0^c \\ v_{ij}^d \geq d_{ij} + T_0^d - \hat{u}_{ij}^d \\ v_{ij}^d \geq \hat{u}_{ij}^d - d_{ij} - T_0^d \\ u_k^a \leq u_k^b \leq u_k^c \leq u_k^d \quad \forall k = \{1,...,n\} \end{cases} \quad \textbf{(L-P}_0\textbf{)}$$

$$\min_{\tilde{u}_1,...,\tilde{u}_n \in \mathcal{T}_n(\mathbb{R}_0^+)} \tfrac{1}{4}\sum_{i=1}^n \sum_{j=1}^n (v_{ij}^a + v_{ij}^b + v_{ij}^c + v_{ij}^d)$$

$$s.t. \begin{cases} v_{ij}^a \geq a_{ij} + T_0^a - \hat{u}_{ij}^a \\ v_{ij}^a \geq \hat{u}_{ij}^a - a_{ij} - T_0^a \\ v_{ij}^b \geq b_{ij} + T_0^b - \hat{u}_{ij}^b \\ v_{ij}^b \geq \hat{u}_{ij}^b - b_{ij} - T_0^b \\ v_{ij}^c \geq c_{ij} + T_0^c - \hat{u}_{ij}^c \\ v_{ij}^c \geq \hat{u}_{ij}^c - c_{ij} - T_0^c \\ v_{ij}^d \geq d_{ij} + T_0^d - \hat{u}_{ij}^d \\ v_{ij}^d \geq \hat{u}_{ij}^d - d_{ij} - T_0^d \\ \sum_{k=1}^n u_k^a = \sigma_a \\ \sum_{k=1}^n u_k^b = \sigma_b \\ \sum_{k=1}^n u_k^c = \sigma_c \\ \sum_{k=1}^n u_k^d = \sigma_d \\ 0 \leq u_k^a \leq u_k^b \leq u_k^c \leq u_k^d \quad \forall k = \{1,...,n\} \end{cases} \quad \textbf{(L-P}_\sigma\textbf{)}$$

**Theorem 3.** For any given TrFPR $\tilde{X} \in \mathbb{A}_n(\tilde{T}_0)$, both the optimization models **P** and **P_σ** (for any predefined total utility constraint $\tilde{\sigma} \in \mathcal{T}(\mathbb{R}^+)$) admit at least one optimal solution $\tilde{u}^*$.

**Proof.** The existence of an optimal solution is guaranteed by the fundamental theorem of linear programming (Dantzig, 1963), provided the model's feasible region is non-empty and its objective function is bounded below.

The feasible region $\Omega$ for **L-P** is: 1) *Non-empty*: A feasible solution can be constructed, for instance, by setting $\tilde{u}_k = T(0,0,0,0)$ for all $k = \{1,...,n\}$, which satisfies all constraints of **L-P**, hence $\Omega \neq \varnothing$. 2) *Bounded*: The objective function $\min \frac{1}{4}\sum_{i=1}^{n}\sum_{j=1}^{n}(v_{ij}^a + v_{ij}^b + v_{ij}^c + v_{ij}^d)$ is bounded below by 0, as all constraints enforce $v_{ij}^a, v_{ij}^b, v_{ij}^c, v_{ij}^d \geq 0$.

The feasible region $\Omega_\sigma$ for **L-P**$_\sigma$ is: 1) *Non-empty*: A feasible solution can be constructed, for instance, by setting $\tilde{u}_k = T(\frac{\sigma_a}{n}, \frac{\sigma_b}{n}, \frac{\sigma_c}{n}, \frac{\sigma_d}{n})$ for all $k = \{1,...,n\}$, which satisfies all constraints of **L-P**$_\sigma$, hence $\Omega_\sigma \neq \varnothing$. 2) *Bounded*: The objective function is the same as that of **P**, which is bounded below by 0. Therefore, both **P** and **P**$_\sigma$ admit at least one optimal solution. $\square$

**Corollary 2.** Let $\tilde{X} \in \mathbb{A}_n(\tilde{T}_0)$ be a consistent TrFPR. Then its LAD utility vector $\tilde{u}$ can be obtained constructively by $\tilde{x}_{ij} \oplus \tilde{T}_0 = \tilde{u}_i \oplus \tilde{u}_j^\circ$ in Theorem 1 without solving the optimization model **P**.

**Proof.** According to Theorem 1, if a TrFPR $\tilde{X}$ is consistent, then there exists a fuzzy utility vector $\tilde{u} \in \mathcal{T}_n([0,1])$ such that $\tilde{x}_{ij} \oplus \tilde{T}_0 = \tilde{u}_i \oplus \tilde{u}_j^\circ, \forall i,j = \{1,...,n\}$. Let us pick $\tilde{u}_i := \tilde{x}_{ik}$ for some $k = \{1,...,n\}$. Then, we can obtain the utility vector $\tilde{u} = (\tilde{x}_{1k},...,\tilde{x}_{nk})$ with the minimum objective function value of 0. Since each $\tilde{x}_{ik} \geq 0$, the constructed utility vector $\tilde{u}$ is a non-negative and feasible optimal solution for **P**, and it achieves the minimum objective value. Therefore, when $\tilde{X}$ is consistent, solving the optimization model (**P**) is unnecessary, as the utility vector $\tilde{u}$ can be obtained simply and directly by $\tilde{x}_{ij} \oplus \tilde{T}_0 = \tilde{u}_i \oplus \tilde{u}_j^\circ$ in Theorem 1. $\square$

**Corollary 3** (*Aggregation of LAD utilities*). Let $\tilde{X}^k = \{\tilde{X}^1,...,\tilde{X}^m\} \in \mathbb{A}_n(\tilde{T}_0)$ be a set of $m$ TrFPRs and $w = \{w_1,...,w_m\}$ be a vector of weights such that $w_k \geq 0$ and $\sum_{k=1}^{m} w_k = 1$. For each $k \in \{1,...,m\}$, let $\tilde{u}^{k*}$ be the optimal LAD utility vector for $\tilde{X}^k$ with the optimal objective value $Z^*(\tilde{X}^k)$. The aggregated utility vector $\tilde{u}^c = \bigoplus_{k=1}^{m} w_k \tilde{u}^{k*}$, associated with the aggregated TrFPR $\tilde{X}^c = \bigoplus_{k=1}^{m} w_k \tilde{X}^k$, yields an objective value $Z(\tilde{X}^c, \tilde{u}^c)$ that satisfies:

i) $\tilde{u}^c$ is a feasible solution for **P** applied to the aggregated TrFPR $\tilde{X}^c$.

ii) The optimal objective value for $\tilde{X}^c$, denoted as $Z^*(\tilde{X}^c)$, is bounded above by the aggregated solution objective value $Z(\tilde{X}^c, \tilde{u}^c)$: $Z^*(\tilde{X}^c) \leq Z(\tilde{X}^c, \tilde{u}^c)$.

iii) The objective value of the aggregated solution $Z(\tilde{X}^c, \tilde{u}^c)$ is bounded above by the weighted sum of individual optimal objective values $\sum_{k=1}^{m} w_k Z^*(\tilde{X}^k)$: $Z(\tilde{X}^c, \tilde{u}^c) \leq \sum_{k=1}^{m} w_k Z^*(\tilde{X}^k)$.

**Proof.** i) To prove that $\tilde{u}^c$ is feasible for **P**, we only need to check the non-negativity constraint. Since each $\tilde{u}^{k*}$ is a solution for **P**, all its components are non-negative. As the weights $w_k$ are also non-negative, the aggregated utility vector $\tilde{u}^c = \bigoplus_{k=1}^{m} w_k \tilde{u}^{k*}$ is a non-negative linear combination of non-negative vectors and is therefore itself non-negative. Thus, $\tilde{u}^c$ is a feasible solution.

ii) Since $\tilde{u}^c$ is feasible, it holds that $Z^*(\tilde{X}^c) \leq Z(\tilde{X}^c, \tilde{u}^c)$.

iii) The objective function is defined as $Z(\tilde{X}, \tilde{u}) = \sum_{i=1}^{n} \sum_{j=1}^{n} d(\tilde{x}_{ij} \oplus \tilde{T}_0, \tilde{u}_i \oplus \tilde{u}_j^\circ)$. To facilitate component-wise analysis, let us index the four components of any TrFN $T(a,b,c,d)$ by $\alpha = \{1,2,3,4\}$, such that $T^1 = a, T^2 = b, T^3 = c, T^4 = d$. Note that $\tilde{X}^c = \bigoplus_{k=1}^{m} w_k \tilde{X}^k$ and $\tilde{u}^c = \bigoplus_{k=1}^{m} w_k \tilde{u}^{k*}$. Then

$$d\left(\tilde{x}_{ij}^c \oplus \tilde{T}_0, \tilde{u}_i^c \oplus (\tilde{u}_j^c)^\circ\right) = \tfrac{1}{4}\sum_{\alpha=1}^{4}\left|(\alpha_{ij}^c + T_0^\alpha) - (u_i^{c,\alpha} + 1 - u_j^{c,5-\alpha})\right|$$

$$= \tfrac{1}{4}\sum_{\alpha=1}^{4}\left|\sum_{k=1}^{m} w_k \alpha_{ij}^k + \sum_{k=1}^{m} w_k T_0^\alpha - \left(\sum_{k=1}^{m} w_k u_i^{k*,\alpha} + \sum_{k=1}^{m} w_k - \sum_{k=1}^{m} w_k u_j^{k*,5-\alpha}\right)\right|$$

$$= \tfrac{1}{4}\sum_{\alpha=1}^{4}\left|\sum_{k=1}^{m} w_k \left(\alpha_{ij}^k + T_0^\alpha - \left(u_i^{k*,\alpha} + 1 - u_j^{k*,5-\alpha}\right)\right)\right|$$

$$\leq \sum_{k=1}^{m} w_k d\left(\tilde{x}_{ij}^k \oplus \tilde{T}_0, \tilde{u}_i^{k*} \oplus (\tilde{u}_j^{k*})^\circ\right).$$

Therefore, $Z(\tilde{X}^c, \tilde{u}^c) \leq \sum_{k=1}^{m} w_k Z^*(\tilde{X}^k)$. $\square$

This result indicates that aggregating the individual LAD utility vectors yields a collective utility vector, such that the LAD involved in approximating the collective opinion is bounded by the LAD associated with the individual TrFPRs. In particular, when the individual TrFPRs are consistent, this aggregation preserves consistency and maintains the exact relationship between aggregated utilities and preferences.

**Proposition 6.** Let $\tilde{X}^k = \{\tilde{X}^1,...,\tilde{X}^m\} \in \mathbb{A}_n(\tilde{T}_0)$ be a set of $m$ consistent TrFPRs and let $\tilde{u}^k$ be the utility vector for each $\tilde{X}^k$, such that $\tilde{x}_{ij}^k \oplus \tilde{T}_0 = \tilde{u}_i^k \oplus (\tilde{u}_j^k)^\circ$. Let $\tilde{X}^c = \bigoplus_{k=1}^m w_k \tilde{X}^k$ be the aggregated TrFPR with weights $w_k \geq 0$ and $\sum_{k=1}^m w_k = 1$. Then, the aggregated matrix $\tilde{X}^c$ is also consistent and the utility vector corresponding to $\tilde{X}^c$, denoted as $\tilde{u}^c$, is equal to the weighted average of the individual utility vectors: $\tilde{u}^c = \bigoplus_{k=1}^m w_k \tilde{u}^{k*}$.

**Proof.** Let $\tilde{u}^c = \bigoplus_{k=1}^m w_k \tilde{u}^{k*}$ denote as the aggregated utility vector. Since each $\tilde{X}^k$ is consistent, it holds that $\tilde{x}_{ij}^k \oplus \tilde{T}_0 = \tilde{u}_i^k \oplus (\tilde{u}_j^k)^\circ$. Then, it follows that

$$\tilde{x}_{ij}^c \oplus \tilde{T}_0 = (\bigoplus_{k=1}^m w_k \tilde{x}_{ij}^k) \oplus \tilde{T}_0 = \bigoplus_{k=1}^m w_k (\tilde{x}_{ij}^k \oplus \tilde{T}_0)$$
$$= \bigoplus_{k=1}^m w_k (\tilde{u}_i^k \oplus (\tilde{u}_j^k)^\circ) = (\bigoplus_{k=1}^m w_k \tilde{u}_i^k) \oplus (\bigoplus_{k=1}^m w_k (\tilde{u}_j^k)^\circ)$$
$$= \tilde{u}_i^c \oplus (\tilde{u}_j^c)^\circ \quad \forall i, j \in \{1,...,n\}, k \in \{1,...,m\}.$$

This is the consistency condition for $\tilde{X}^c$ with the utility vector $\tilde{u}^c$. This proves both that $\tilde{X}^c$ is consistent and its corresponding utility vector is $\tilde{u}^c = \bigoplus_{k=1}^m w_k \tilde{u}^{k*}$. □

Corollary 3 guarantees that, in group decision-making (GDM) problems (Li et al., 2023), the collective LAD utility vector can be computed by aggregating the individual LAD utility vectors, which is a critical step in the consensus-reaching process (Lu et al., 2023). Furthermore, Proposition 6 proves that under a consistent TrFPR, this aggregation could preserve the consistency and structure of the preferences. These findings guarantee that aggregating individual LAD utility vectors is a robust and theoretically well-founded method for deriving the collective utility in GDM problems.

### 3.3. Least absolute deviations utility for TrMPRs

While the previous sections have developed a comprehensive LAD utility model for TrFPRs, which are defined on an additive scale, DMs often express their preferences on a multiplicative scale (e.g., Saaty's 1-9 scale). To enhance the applicability of the LAD utility model, this section extends it to TrMPRs using a bijective mapping $\Phi : \mathbb{A}_n(\tilde{T}_0) \to \mathbb{M}_n^m(\tilde{S}_0)$, which is introduced in Proposition 1. To support the definition of TrMPRs, we introduce the following set of operators.

**Definition 14** (*Multiplication*). Let $\mathcal{T}(\mathbb{R}^+)$ be a set of TrFNs defined on the set of positive real numbers. For any two TrFNs $\tilde{T}_1 \equiv T(a_1,b_1,c_1,d_1) \in \mathcal{T}(\mathbb{R}^+)$ and $\tilde{T}_2 \equiv T(a_2,b_2,c_2,d_2) \in \mathcal{T}(\mathbb{R}^+)$, their product $*$ is defined as $\tilde{T}_1 * \tilde{T}_2 = (a_1 a_2, b_1 b_2, c_1 c_2, d_1 d_2)$.

**Proposition 7.** The set of TrFNs on the positive real line, $\mathcal{T}(\mathbb{R}^+)$, together with the multiplicative operation $*$ defined in Definition 14, forms a monoid, denoted as $(\mathcal{T}(\mathbb{R}^+), *)$.

**Proof.** Let $\tilde{T}_1 \equiv T(a_1,b_1,c_1,d_1) \in \mathcal{T}(\mathbb{R}^+)$, $\tilde{T}_2 \equiv T(a_2,b_2,c_2,d_2) \in \mathcal{T}(\mathbb{R}^+)$ and $\tilde{T}_3 \equiv T(a_3,b_3,c_3,d_3) \in \mathcal{T}(\mathbb{R}^+)$ be arbitrary TrFNs in $\mathcal{T}(\mathbb{R}^+)$. The $(\mathcal{T}(\mathbb{R}^+), *)$ satisfies: i) Closure: $\forall \tilde{T}_1, \tilde{T}_2 \in \mathcal{T}(\mathbb{R}^+)$, $\tilde{T}_1 * \tilde{T}_2 \in \mathcal{T}(\mathbb{R}^+)$. ii) Associativity: $(\tilde{T}_1 * \tilde{T}_2) * \tilde{T}_3 = T(a_1 a_2, b_1 b_2, c_1 c_2, d_1 d_2) * T(a_3,b_3,c_3,d_3) = T(a_1 a_2 a_3, b_1 b_2 b_3, c_1 c_2 c_3, d_1 d_2 d_3)$
$= T(a_1,b_1,c_1,d_1) * T(a_2 a_3, b_2 b_3, c_2 c_3, d_2 d_3) = \tilde{T}_1 * (\tilde{T}_2 * \tilde{T}_3)$. iii) Identity element: Let us consider $\tilde{T}_e = T(1,1,1,1) \in \mathcal{T}(\mathbb{R}^+)$, such that $\tilde{T}_1 * \tilde{T}_e = T(a_1,b_1,c_1,d_1) * T(1,1,1,1) = T(a_1,b_1,c_1,d_1) = \tilde{T}_1$. Since the $(\mathcal{T}(\mathbb{R}^+), *)$ is close, associative and has an identity element, $(\mathcal{T}(\mathbb{R}^+), *)$ is a monoid. □

**Definition 15** (*Inverse for TrFNs*). Let $\tilde{T} \equiv T(a,b,c,d) \in \mathcal{T}(\mathbb{R}^+)$ be a TrFN. The inverse function $I : \mathcal{T}(\mathbb{R}^+) \to \mathcal{T}(\mathbb{R}^+)$ of the TrFN $\tilde{T}$ is defined as:

$$I(\tilde{T}) = I(T(a,b,c,d)) = T(\tfrac{1}{d},\tfrac{1}{c},\tfrac{1}{b},\tfrac{1}{a}).$$

**Proposition 8** (*Involutive Property*). The multiplicative inverse operator $I(\cdot)$ is involutive. For any $\tilde{T} \in \mathcal{T}(\mathbb{R}^+)$, it holds that $I(I(\tilde{T})) = \tilde{T}$.

**Proof.** Let $\tilde{T} \equiv T(a,b,c,d)$, then we obtain: $I(I(\tilde{T})) = I(T(\tfrac{1}{d},\tfrac{1}{c},\tfrac{1}{b},\tfrac{1}{a})) = T(a,b,c,d) = \tilde{T}$. □

While the inverse $I(\cdot)$ and negation $N_S(\cdot)$ operators belong to different scale systems, they are not independent. The mapping function $\phi$ establishes a formal isomorphism between them, ensuring that the structure is preserved across the two domains. This preservation is formally stated as follows.

**Proposition 9.** For any TrFN $\tilde{x} \in \mathcal{T}([0,1])$ and its standard negation $N_S(\tilde{x})$, the following holds: $\phi(N_S(\tilde{x})) = I(\phi(\tilde{x}))$.

**Proof.** Let $\tilde{x} = T(a,b,c,d)$. By definition 8, its standard negation $N_S(\tilde{x})$ is given by $T(1-d, 1-c, 1-b, 1-a)$, then

$$\begin{aligned}
\phi(N_S(\tilde{x})) &= T(\phi(1-d), \phi(1-c), \phi(1-b), \phi(1-a)) \\
&= T(m^{2(1-d)-1}, m^{2(1-c)-1}, m^{2(1-b)-1}, m^{2(1-a)-1}) \\
&= T(m^{-(2d-1)}, m^{-(2c-1)}, m^{-(2b-1)}, m^{-(2a-1)}) \\
&= T(\tfrac{1}{\phi(d)}, \tfrac{1}{\phi(c)}, \tfrac{1}{\phi(b)}, \tfrac{1}{\phi(a)}) = I(\phi(\tilde{x})).
\end{aligned}$$

□

The mapping $\phi$ not only links the negation and inverse operators but also preserves the property of neutrality. The additive neutral preference $\tilde{T}_0$ is defined by the fixed-point property $\tilde{T}_0^\circ = \tilde{T}_0$. Applying the identity from Proposition 9, $\phi(\tilde{T}^\circ) = I(\phi(\tilde{T}))$, to this property yields: $I(\phi(\tilde{T}_0)) = \phi(\tilde{T}_0^\circ) = \phi(\tilde{T}_0)$.

This result indicates that $\phi(\tilde{T}_0)$ is a fixed point of the multiplicative inverse operator $I$, which can be used to define the neutral TrFN in the multiplicative framework.

**Definition 16** (*Multiplicative neutral TrFN*). A TrFN $\tilde{S}_0 \in \mathcal{T}(\mathbb{R}^+)$ is called a multiplicative neutral TrFN if it satisfies the fixed-point property with respect to the inverse operator:
$$I(\tilde{S}_0) = \tilde{S}_0.$$
A multiplicative neutral TrFN may be obtained by transforming an additive neutral TrFN: $\tilde{S}_0 = \phi(\tilde{T}_0)$.

It is critical to note that since $\tilde{T}_0$ is an element of $\mathcal{T}([0,1])$, its image $\tilde{S}_0$ under the mapping $\phi$ will be an element of $\mathcal{T}([\tfrac{1}{m}, m])$. This means the multiplicative neutral preference, like the preferences in Saaty's classical scale, is inherently bounded.

**Example 7.** Given the additive neutral TrFN $\tilde{T}_0 = T(0.4, 0.5, 0.5, 0.6)$ from Example 3 and setting the multiplicative scale $m=9$, we can compute its multiplicative neutral TrFN by the bijective mapping $\phi$ in Proposition 1:
$$\tilde{S}_0 = \phi(\tilde{T}_0) = T\left(9^{2\times 0.4-1}, 9^{2\times 0.5-1}, 9^{2\times 0.5-1}, 9^{2\times 0.6-1}\right) = T(9^{-0.2}, 1, 1, 9^{0.2}).$$

This satisfies the fixed-point property (Definition 16): $I(\tilde{S}_0) = I(T(9^{-0.2}, 1, 1, 9^{0.2})) = T(\tfrac{1}{9^{0.2}}, \tfrac{1}{1}, \tfrac{1}{1}, \tfrac{1}{9^{-0.2}})$ $= T(9^{-0.2}, 1, 1, 9^{0.2}) = \tilde{S}_0$. Based on the multiplicative neutral preference $\tilde{S}_0$ and the inverse function $I$, we propose the definition of TrMPR.

**Definition 17** (*Multiplicative trapezoidal fuzzy preference relation (TrMPR)*). Let $n \in \mathbb{N}$. A preference matrix $\tilde{Y} = (\tilde{y}_{ij}) \in \mathcal{M}_{n\times n}(\mathcal{T}), \forall i,j = \{1,\ldots,n\}$, where $\tilde{y}_{ij} = (y_{ij}^a, y_{ij}^b, y_{ij}^c, y_{ij}^d) \in \mathcal{T}(\mathbb{R}^+)$ is a TrMPR if

(i) $\tilde{y}_{kk} = \tilde{S}_0 \ \forall k \in \{1,\ldots,n\}$,

(ii) $\tilde{y}_{ij} = I(\tilde{y}_{ji}) \ \forall i,j \in \{1,\ldots,n\}, i \neq j$,

where $\tilde{S}_0 \in \mathcal{T}(\mathbb{R}^+)$ represents the multiplicative neutral preference and satisfies $I(\tilde{S}_0) = \tilde{S}_0$. If $\tilde{y}_{ij} \in \left[\frac{1}{m}, m\right], \forall i,j \in \{1,\ldots,n\}$, where $m \in \mathbb{N} \setminus \{1\}$, we say that the TrMPR has been elicited in the scale $m$. Usually, $m=9$, but it could be any other natural number greater than one. The set of TrMPRs of dimension $n$ under scale $m$ with multiplicative neutral preference $\tilde{S}_0$ is denoted as $\mathbb{M}_n^m(\tilde{S}_0)$.

**Example 8.** Let us transform the TrFPR $\tilde{X}$ from Example 3 into its corresponding TrMPR $\tilde{Y} = \Phi(\tilde{X})$ with $m=9$:

$$\tilde{Y} = \Phi(\tilde{X}) = \begin{pmatrix} \phi(\tilde{x}_{11}) & \phi(\tilde{x}_{12}) & \phi(\tilde{x}_{13}) \\ \phi(\tilde{x}_{21}) & \phi(\tilde{x}_{22}) & \phi(\tilde{x}_{23}) \\ \phi(\tilde{x}_{31}) & \phi(\tilde{x}_{32}) & \phi(\tilde{x}_{33}) \end{pmatrix} = \begin{pmatrix} T(9^{-0.2},1,1,9^{0.2}) & T(9^{0.2},9^{0.4},9^{0.4},9^{0.6}) & T(9^{0.2},9^{0.4},9^{0.6},9^{0.8}) \\ T(9^{-0.6},9^{-0.4},9^{-0.4},9^{-0.2}) & T(9^{-0.2},1,1,9^{0.2}) & T(1,9^{0.2},9^{0.4},9^{0.6}) \\ T(9^{-0.8},9^{-0.6},9^{-0.4},9^{-0.2}) & T(9^{-0.6},9^{-0.4},9^{-0.2},1) & T(9^{-0.2},1,1,9^{0.2}) \end{pmatrix}.$$

This matrix satisfies: i) $\tilde{y}_{11} = \tilde{y}_{22} = \tilde{y}_{33} = T(9^{-0.2},1,1,9^{0.2}) = \tilde{S}_0$. ii) $I(\tilde{y}_{12}) = I(T(9^{0.2},9^{0.4},9^{0.6},9^{0.8})) = T(9^{-0.8},9^{-0.6},9^{-0.4},9^{-0.2}) = \tilde{y}_{21}$. Similarly, $I(\tilde{y}_{13}) = \tilde{y}_{31}$ and $I(\tilde{y}_{23}) = \tilde{y}_{32}$.

**Definition 18** (*Consistent TrMPR*). A TrMPR $\tilde{Y} = (\tilde{y}_{ij}) \in \mathbb{M}_n^m(\tilde{S}_0)$ is said to be consistent if it satisfies

$$\tilde{y}_{ij} * \tilde{S}_0 = \tilde{y}_{ik} * \tilde{y}_{kj} \ \forall i,j,k \in \{1,\ldots,n\}$$

where the operation $*$ is defined in Definition 14. When the neutral reference TrFN $\tilde{T}_0 = T(\frac{1}{2},\frac{1}{2},\frac{1}{2},\frac{1}{2})$, the multiplicative neutral $\tilde{S}_0 = \phi(\tilde{T}_0)$ becomes: $\tilde{S}_0 = \phi(T(\frac{1}{2},\frac{1}{2},\frac{1}{2},\frac{1}{2})) = T(m^{2\cdot\frac{1}{2}-1}, m^{2\cdot\frac{1}{2}-1}, m^{2\cdot\frac{1}{2}-1}, m^{2\cdot\frac{1}{2}-1}) = T(1,1,1,1)$, which corresponds to the crisp value 1. If all input preference values $\tilde{y}_{ij}$ are crisp numbers, then both the definition of TrMPR and consistent TrMPR reduce to the classical definitions of MPR (Saaty, 2008) and consistent MPR (Saaty, 2008).

**Example 9.** Consider the TrMPR $\tilde{Y} = (\tilde{y}_{ij})$ from Example 8. We check if it satisfies the consistency condition: $\tilde{y}_{ij} * \tilde{S}_0 = \tilde{y}_{ik} * \tilde{y}_{kj}, \forall i,j,k \in \{1,\ldots,n\}$. Specifically, for $i=1, j=2, k=3$:

$$\tilde{y}_{12} * \tilde{S}_0 = T(9^{0.2},9^{0.4},9^{0.4},9^{0.6}) * T(9^{-0.2},1,1,9^{0.2}) = T(1,9^{0.4},9^{0.6},9),$$

$$\tilde{y}_{13} * \tilde{y}_{32} = T(9^{0.2},9^{0.4},9^{0.6},9^{0.8}) * T(9^{-0.6},9^{-0.4},9^{-0.2},1) = T(9^{-0.4},1,9^{0.4},9^{0.8}).$$

Since $\tilde{y}_{ij} * \tilde{S}_0 \neq \tilde{y}_{ik} * \tilde{y}_{kj}$, we conclude that the TrMPR $\tilde{Y}$ is inconsistent same as the TrFPR $\tilde{X}$ in Example 4.

**Theorem 4.** Let $n \in \mathbb{N}$ and consider $m \in \mathbb{N} \setminus \{1\}$. A TrFPR $\tilde{X} \in \mathbb{A}_n(\tilde{T}_0)$ is consistent if and only if the TrMPR $\Phi(\tilde{X}) \in \mathbb{M}_n^m(\tilde{S}_0)$ is consistent.

**Proof.** By definition 12, $\tilde{X}$ is consistent if $\tilde{x}_{ij} \oplus \tilde{T}_0 = \tilde{x}_{ik} \oplus \tilde{x}_{kj}$, for any $i,j,k \in \{1,\ldots,n\}$, the following holds:

$$\tilde{x}_{ij} \oplus \tilde{T}_0 = \tilde{x}_{ik} \oplus \tilde{x}_{kj} \ \forall i,j,k \in \{1,\ldots,n\}$$
$$\Leftrightarrow \alpha_{ij} + T_0^\alpha = \alpha_{ik} + \alpha_{kj} \ \forall \alpha = \{a,b,c,d\}, i,j,k \in \{1,\ldots,n\}$$
$$\Leftrightarrow 2\alpha_{ij} + 2T_0^\alpha - 2 = 2\alpha_{ik} + 2\alpha_{kj} - 2 \ \forall \alpha = \{a,b,c,d\}, i,j,k \in \{1,\ldots,n\}$$
$$\Leftrightarrow m^{(2\alpha_{ij}-1)+(2T_0^\alpha-1)} = m^{(2\alpha_{ik}-1)+(2\alpha_{kj}-1)} \ \forall \alpha = \{a,b,c,d\}, i,j,k \in \{1,\ldots,n\}$$
$$\Leftrightarrow m^{(2\alpha_{ij}-1)} \cdot m^{(2T_0^\alpha-1)} = m^{(2\alpha_{ik}-1)} \cdot m^{(2\alpha_{kj}-1)} \ \forall \alpha = \{a,b,c,d\}, i,j,k \in \{1,\ldots,n\}$$
$$\Leftrightarrow \tilde{y}_{ij} * \tilde{S}_0 = \tilde{y}_{ik} * \tilde{y}_{kj} \ \forall i,j,k \in \{1,\ldots,n\}.$$

which is the multiplicative consistency for $\Phi(\tilde{X})$. □

The previous result allow us to extend the LAD utility vector concept to TrMPRs for pairwise comparisons given in a multiplicative scale. Assume $\tilde{Y} \in \mathbb{M}_n^m(\tilde{S}_0)$ is a TrMPR for some $m \in \mathbb{N} \setminus \{1\}$. Then, we can consider the TrFPR $\tilde{X} = \Phi^{-1}(\tilde{Y})$, where $\Phi^{-1}: \mathbb{M}_n^m(\tilde{S}_0) \to \mathbb{A}_n(\tilde{T}_0)$, and compute the corresponding LAD utility vector $\tilde{u} \in \mathcal{T}_n(\mathbb{R}^+)$. Therefore, we define the multiplicative LAD utility vector as follows.

**Definition 19** (*LAD utility for TrMPRs*). Let $\tilde{Y} = (\tilde{y}_{ij}) \in \mathbb{M}_n^m(\tilde{S}_0)$ be a TrMPR under the scale $m \in \mathbb{N} \setminus \{1\}$. An LAD utility vector for $\tilde{Y}$ is a trapezoidal fuzzy vector $\tilde{u} = (\tilde{u}_1, \tilde{u}_2, ..., \tilde{u}_n) \in \mathcal{T}_n(\mathbb{R}^+)$ that minimizes the total absolute deviation between the shifted TrFPR $\tilde{x}_{ij} \oplus \tilde{T}_0$ transformed from the original TrMPR and that generated by the utility vector, i.e.,

$$\min_{\tilde{u}_1,...,\tilde{u}_n \in \mathcal{T}_n(\mathbb{R}^+)} \sum_{i=1}^n \sum_{j=1}^n d\left(\phi^{-1}(\tilde{y}_{ij}) \oplus \phi^{-1}(\tilde{S}_0), \tilde{u}_i \oplus \tilde{u}_j^\circ\right)$$

subject to:

$$0 \leq u_k^a \leq u_k^b \leq u_k^c \leq u_k^d \quad \forall k = \{1,...,n\},$$

where $d(\cdot, \cdot)$ is the normalized Manhattan (L1) distance for TrFNs (Definition 3), $\tilde{u}_j^\circ$ is the negation of $\tilde{u}_j$, and $\tilde{S}_0$ is the reference TrFN representing the multiplicative neutral preference.

The LAD utility vector $\tilde{u}$ provides an absolute measure of preference intensity for each alternative. However, in many decision-making contexts such as the AHP, it is common to express these preferences as a weight vector $w$, where the sum of the weights equals one in crisp terms. Although the literature on fuzzy weights has proposed various normalization approaches, there is still no widely accepted and clearly defined method (Wang and Elhag, 2006). It is worth noting that under fuzzy environments, the sum of fuzzy weights cannot always be strictly equal to one. To inherit the idea of the sum of weights being 1 from the classical definition of weights while maintaining the structural characteristics of TrFN, we introduce a total utility constraint $\tilde{\sigma}$. Then the utility vector $\tilde{u} = (\tilde{u}_1, \tilde{u}_2, ..., \tilde{u}_n) \in \mathcal{T}_n(\mathbb{R}^+)$ for a TrMPR can be obtained by solving an optimization model $\mathbf{Q}_\sigma$ under the constraint $\tilde{\sigma}$:

$$\min_{\tilde{u}_1,...,\tilde{u}_n \in \mathcal{T}_n(\mathbb{R}^+)} \sum_{i=1}^n \sum_{j=1}^n d\left(\phi^{-1}(\tilde{y}_{ij}) \oplus \phi^{-1}(\tilde{S}_0), \tilde{u}_i \oplus \tilde{u}_j^\circ\right)$$

$$s.t. \begin{cases} \bigoplus_{k=1}^n \tilde{u}_k = \tilde{\sigma} \\ 0 \leq u_k^a \leq u_k^b \leq u_k^c \leq u_k^d \quad \forall k = \{1,...,n\} \end{cases} \quad (\mathbf{Q}_\sigma)$$

where $\tilde{\sigma}$ is a TrFN approximately equal $T(1,1,1,1)$, which can be subjectively specified by the DM, for example $\tilde{\sigma} = T(0.8, 0.9, 1.1, 1.2)$. Consequently, we formally define the LAD fuzzy weight based on model $\mathbf{Q}_\sigma$.

**Definition 20** (*LAD Fuzzy weight for TrMPR*). For a given TrMPR $\tilde{Y}$, the LAD fuzzy weight vector $\tilde{w} = (\tilde{w}_1, ..., \tilde{w}_n)$ is obtained as the optimal solution to the LAD utility with the total utility constraint $\tilde{\sigma}$, which approximately equal $T(1,1,1,1)$, i.e.,

$$\min_{\tilde{w}_1,...,\tilde{w}_n \in \mathcal{T}_n(\mathbb{R}^+)} \sum_{i=1}^n \sum_{j=1}^n d\left(\phi^{-1}(\tilde{y}_{ij}) \oplus \phi^{-1}(\tilde{S}_0), \tilde{w}_i \oplus \tilde{w}_j^\circ\right)$$

subject to:

$$\bigoplus_{k=1}^n \tilde{w}_k = \tilde{\sigma}, \quad 0 \leq w_k^a \leq w_k^b \leq w_k^c \leq w_k^d \quad \forall k = \{1,...,n\},$$

where $d(\cdot, \cdot)$ is the normalized Manhattan (L1) distance for TrFNs (Definition 3), $\tilde{w}_j^\circ$ is the negation of $\tilde{w}_j$, and $\tilde{S}_0$ is the reference TrFN representing the multiplicative neutral preference. The normalized weight vector $\tilde{w}$ can be used to rank alternatives by computing the magnitude $Mag(\tilde{w})$.

**Example 10.** Consider the TrMPR $\tilde{Y} = (\tilde{y}_{ij})$ from Example 8. First, we obtain the LAD utility

vector by Definition 19, $\tilde{u}^* = (T(0.3,0.3,0.3,0.5), T(0.1,0.1,0.1,0.3), T(0,0,0,0.2))$ with an objective value of 0.2. These results are identical to those in Example 6, which provides empirical support for the coherence of the $\Phi$ mapping. Second, we set the total utility constraint $\tilde{\sigma} = (0.8, 0.9, 1.1, 1.2)$ and derive the LAD fuzzy weight by Definition 20. The resulting optimal LAD fuzzy weight vector is $\tilde{w}^* = (T(0.434, 0.434, 0.534, 0.534), T(0.233, 0.333, 0.333, 0.333), T(0.133, 0.133, 0.233, 0.333))$ with an objective value of 0.6. The results provide the ranking $A_1 \succ A_2 \succ A_3$, the same as Example 6.

## 4. Application and comparative analysis in Fuzzy AHP

This study develops a systematic theoretical framework for deriving the LAD utility vector from a TrFPR. The purpose of this section is to demonstrate the applicability and operational procedure of the proposed model through a classical multi-criteria decision-making (MCDM) scenario using the Fuzzy AHP approach, and to compare it with the traditional Fuzzy AHP method.

### 4.1. Illustrative example: evaluation of land development projects

To demonstrate the applicability and procedure of the proposed LAD utility model, we consider a realistic land development project evaluation problem concerning the marketization of collectively-owned commercial construction land (COCCL). As a crucial land system reform in China, marketization of COCCL breaks the state's monopoly on land development rights by allowing collectively-owned land to enter the market on equal terms with state-owned land. A government aims to evaluate and prioritize four COCCL projects for strategic support and resource allocation:

- Alternative 1: Homestead conversion for rural revitalization ($A_1$),
- Alternative 2: Agri-tourism industry cluster ($A_2$),
- Alternative 3: Mining industry supports infrastructure ($A_3$),
- Alternative 4: Integrated eco-wellness industrial park ($A_4$).

The evaluation is based on three key criteria: economic benefits ($C_1$), social impact ($C_2$), and environmental sustainability ($C_3$). For simplicity, we assume that the criterion weights are given as crisp numbers: $\Omega_1 = 0.5$, $\Omega_2 = 0.3$, $\Omega_3 = 0.2$. A DM provides their judgments by comparing the four alternatives under each criterion, resulting in three TrMPRs $\tilde{Y}^k (k=1,2,3)$. These judgments are expressed on Saaty's 1-9 scale, with $m=9$. We assume a personalized additive neutral preference $\tilde{T}_0 = T(0.4, 0.5, 0.5, 0.6)$ with multiplicative neutral TrFN $\tilde{S}_0 = \phi(\tilde{T}_0) = T(9^{-0.2}, 1, 1, 9^{0.2})$.

$$\tilde{Y}_1 = \begin{pmatrix} T(9^{-0.2},1,1,9^{0.2}) & T(\tfrac{1}{3},\tfrac{1}{2},\tfrac{1}{2},1) & T(2,3,3,4) & T(3,4,5,6) \\ T(1,2,2,3) & T(9^{-0.2},1,1,9^{0.2}) & T(4,5,6,7) & T(5,6,7,8) \\ T(\tfrac{1}{4},\tfrac{1}{3},\tfrac{1}{3},\tfrac{1}{2}) & T(\tfrac{1}{7},\tfrac{1}{6},\tfrac{1}{5},\tfrac{1}{4}) & T(9^{-0.2},1,1,9^{0.2}) & T(1,1,2,3) \\ T(\tfrac{1}{6},\tfrac{1}{5},\tfrac{1}{4},\tfrac{1}{3}) & T(\tfrac{1}{8},\tfrac{1}{7},\tfrac{1}{6},\tfrac{1}{5}) & T(\tfrac{1}{3},\tfrac{1}{2},1,1) & T(9^{-0.2},1,1,9^{0.2}) \end{pmatrix}.$$

$$\tilde{Y}_2 = \begin{pmatrix} T(9^{-0.2},1,1,9^{0.2}) & T(2,3,4,5) & T(5,6,7,8) & T(1,1,2,3) \\ T(\tfrac{1}{5},\tfrac{1}{4},\tfrac{1}{3},\tfrac{1}{2}) & T(9^{-0.2},1,1,9^{0.2}) & T(3,4,4,5) & T(\tfrac{1}{3},\tfrac{1}{2},\tfrac{1}{2},1) \\ T(\tfrac{1}{8},\tfrac{1}{7},\tfrac{1}{6},\tfrac{1}{5}) & T(\tfrac{1}{5},\tfrac{1}{4},\tfrac{1}{4},\tfrac{1}{3}) & T(9^{-0.2},1,1,9^{0.2}) & T(\tfrac{1}{6},\tfrac{1}{5},\tfrac{1}{4},\tfrac{1}{3}) \\ T(\tfrac{1}{3},\tfrac{1}{2},1,1) & T(1,2,2,3) & T(3,4,5,6) & T(9^{-0.2},1,1,9^{0.2}) \end{pmatrix}.$$

$$\tilde{Y}_3 = \begin{pmatrix} T(9^{-0.2},1,1,9^{0.2}) & T(2,3,3,4) & T(6,7,7,8) & T(4,5,5,6) \\ T(\tfrac{1}{4},\tfrac{1}{3},\tfrac{1}{3},\tfrac{1}{2}) & T(9^{-0.2},1,1,9^{0.2}) & T(4,5,5,6) & T(2,3,3,4) \\ T(\tfrac{1}{8},\tfrac{1}{7},\tfrac{1}{7},\tfrac{1}{6}) & T(\tfrac{1}{6},\tfrac{1}{5},\tfrac{1}{5},\tfrac{1}{4}) & T(9^{-0.2},1,1,9^{0.2}) & T(\tfrac{1}{3},\tfrac{1}{3},\tfrac{1}{2},\tfrac{1}{2}) \\ T(\tfrac{1}{6},\tfrac{1}{5},\tfrac{1}{5},\tfrac{1}{4}) & T(\tfrac{1}{4},\tfrac{1}{3},\tfrac{1}{3},\tfrac{1}{2}) & T(2,2,3,3) & T(9^{-0.2},1,1,9^{0.2}) \end{pmatrix}.$$

For each TrMPR $\tilde{Y}^k (k=1,2,3)$, let us consider $\tilde{\sigma} = T(0.8, 0.9, 1.1, 1.2)$ and derive the local fuzzy weights $\tilde{w}_i^k (i=1,2,3,4, k=1,2,3)$ under each criterion by solving the LAD model (Definition 20) using CPLEX 22.1.0 through Python. For each criterion, we obtain

$$\tilde{w}_i^1 = \begin{pmatrix} T(0.3000, 0.3291, 0.3713, 0.3713), T(0.4577, 0.4868, 0.4868, 0.5073), \\ T(0.0423, 0.0790, 0.1628, 0.2423), T(0.0000, 0.0051, 0.0791, 0.0791). \end{pmatrix},$$

$$\tilde{w}_i^2 = \begin{pmatrix} T(0.3711, 0.3711, 0.5155, 0.5712), T(0.1500, 0.2134, 0.2134, 0.2577), \\ T(0.0000, 0.0000, 0.0000, 0.0000), T(0.2789, 0.3155, 0.3711, 0.3711). \end{pmatrix},$$

$$\tilde{w}_i^3 = \begin{pmatrix} T(0.4244, 0.4996, 0.4996, 0.5732), T(0.2667, 0.2667, 0.3667, 0.3667), \\ T(0.0000, 0.0004, 0.0248, 0.0512), T(0.1089, 0.1333, 0.2089, 0.2089). \end{pmatrix}.$$

The global weights for each alternative $A_i$ are obtained by aggregating the local LAD weights using the given criteria weights: $\tilde{w}_i = (0.5\tilde{w}_i^1) \oplus (0.3\tilde{w}_i^2) \oplus (0.2\tilde{w}_i^3)$:

$$\tilde{w}_i = \begin{pmatrix} T(0.3462, 0.3758, 0.4402, 0.4716), T(0.3272, 0.3607, 0.3808, 0.4043), \\ T(0.0212, 0.0396, 0.0834, 0.1314), T(0.1054, 0.1239, 0.1926, 0.1927). \end{pmatrix}.$$

The alternatives are then ranked by applying the $Mag(\cdot)$ function to their final weights as follows:

$$Mag(\tilde{w}^1) = 0.4082, \ Mag(\tilde{w}^2) = 0.3699, \ Mag(\tilde{w}^3) = 0.0652, \ Mag(\tilde{w}^4) = 0.1567.$$

This leads to the final rank order: $A_1 \succ A_2 \succ A_4 \succ A_3$. The proposed LAD model has thus successfully distilled the complex and inconsistent fuzzy judgments into a clear and interpretable ranking, revealing a sophisticated preference structure that provides decision support for the government.

### 4.2. Comparative analysis

This section compares the use of the LAD utility vector for TrMPRs (Definition 20) with the fuzzy arithmetic mean method (AMM) (van Laarhoven and Pedrycz, 1983) and fuzzy geometric mean method (GMM) (Buckley, 1985) to compute priorities from an TrMPR, which is widely applied in fuzzy AHP. We still consider the TrMPRs in Sec. 4.1, taking $\tilde{Y}_1$ as example, the corresponding fuzzy weights by AMM and GMM we compute as in Table 1.

Table 1. The fuzzy weight vectors obtained by various methods.

| Method | Fuzzy weight | Deviation |
| --- | --- | --- |
| LAD (Proposed) | $\begin{pmatrix} T(0.3000, 0.3291, 0.3713, 0.3713), T(0.4577, 0.4868, 0.4868, 0.5073), \\ T(0.0423, 0.0790, 0.1628, 0.2423), T(0.0000, 0.0051, 0.0791, 0.0791). \end{pmatrix}$ | 0.9823 |
| AMM | $\begin{pmatrix} T(0.1476, 0.2703, 0.3539, 0.6298), T(0.2629, 0.4452, 0.5961, 0.9811), \\ T(0.0503, 0.0795, 0.1316, 0.2660), T(0.0316, 0.0586, 0.0900, 0.1548). \end{pmatrix}$ | 2.3955 |
| GMM | $\begin{pmatrix} T(0.1344, 0.2727, 0.3195, 0.6846), T(0.2390, 0.4850, 0.5845, 1.1136), \\ T(0.0491, 0.0846, 0.1167, 0.2421), T(0.0327, 0.0602, 0.0873, 0.1572). \end{pmatrix}$ | 2.6843 |

Table 1 shows that the fuzzy weight vectors derived by GMM could contain components with upper bounds greater than one (e.g., $\tilde{w}_2^1$ from GMM has an upper bound of 1.1136). While mathematically valid, such weights can be less intuitive to interpret, as weights are often expected to fall within the [0,1] interval. The result not only compares the weight vectors obtained by the three different methods but also verifies that the LAD weighting vector proposed in this study is the closest one to the original matrix $\tilde{Y}_1$. Therefore, our methodology prioritizes the preservation of the DM's original judgments, providing a more robust and reliable foundation for subsequent decision-making.

### 5. Conclusion

This study has proposed a novel definition of TrFPR that relies on a neutral TrFN to capture

DMs' perceptions of fuzzy preferences in a more flexible way. We have also developed a method to derive utility vectors from TrFPRs based on the LAD. Such a method is based on an optimization model that finds a utility vector that minimizes the LAD between the shifted original TrFPR and that generated by the utility vector. The discussion highlighted several significant properties and benefits of this approach, which collectively enhance its practical application in decision-making processes

The theoretical contributions of this study are threefold. First, we moved beyond the conventional crisp neutral values by introducing the concept of a neutral TrFN $\tilde{T}_0$. This not only allows DMs to express indifference as a TrFN that better reflects human cognitive reality, but also generalizes the definition of TrFPR. Second, we proposed a theorem for directly extracting utility vectors from consistent TrFPRs, as well as an LAD model for deriving utility vectors from inconsistent TrFPRs. Third, we extended the approach beyond additive TrFPRs by introducing a mapping function and proving that the LAD model can be applied to TrMPRs, thereby establishing a unified theoretical framework. This allows any TrFN to be used as input in the TrFPR and also permits any symmetric TrFN to represent neutrality.

The LAD method streamlines the computational process by reducing the dimensionality of preference modeling from $n^2$ pairwise comparisons to an $n$-dimensional utility vector for each DM. This drastic dimensionality reduction may be particularly advantageous in large-scale group decision-making contexts, as it replaces the complex aggregation of large fuzzy matrices with the far more efficient aggregation of simple utility vectors. Therefore, in the future, we will apply the proposed approach to address large-scale problems, exploring how to effectively aggregate individual LAD utility vectors and manage consensus-reaching processes. Additionally, future work could also focus on developing guidelines for selecting an appropriate $\tilde{\sigma}$ of model $\mathbf{P}_\sigma$ in different decision contexts. Exploring these extensions would significantly enhance the practical scope and theoretical depth of the LAD utility model.

**Acknowledgements**

This research was supported by the Social Science Foundation of Shaanxi Province (No. 2023R051), the Xi'an Association for Science and Technology (No. 25JCZX015R2) and the Xi'an Municipal Social Science Planning Fund (No. 25JX200). This research was also partially supported by the projects ProyExcel_00257 of the Regional Ministry of Economy, Innovation, Science and Employment of Andalusia and PID2024-161073NB-I00 of the Spanish Ministry of Economy and Finance and ERDF.